\input amstex\documentstyle {amsppt}  
\pagewidth{12.5 cm}\pageheight{19 cm}\magnification\magstep1
\topmatter
\title Parabolic character sheaves, I\endtitle
\author G. Lusztig\endauthor
\address Department of Mathematics, M.I.T., Cambridge, MA 02139\endaddress
\dedicatory{Dedicated to Pierre Cartier on the occasion of his 70th birthday}
\enddedicatory
\thanks Supported in part by the National Science Foundation\endthanks
\endtopmatter   
\document 
\define\dsv{\dashv}

\define\po{\text{\rm pos}}

\define\si{\sim}
\define\wt{\widetilde}
\define\sqc{\sqcup}

\define\qua{\quad}

\define\hG{\hat G}
\define\hX{\hat X}
\define\hY{\hat Y}
\define\hcl{\hat{\cl}}
\define\dx{\dot x}
\define\tcl{\ti\cl}

\define\bK{\bar K}
\define\bY{\bar Y}

\define\bvt{\bar\vt}
\define\bpi{\bar\p}
\define\lb{\linebreak}

\define\op{\oplus}
   
\redefine\sp{\spadesuit}
\define\em{\emptyset}
\define\imp{\implies}
\define\ra{\rangle}

\define\iy{\infty}
\define\m{\mapsto}
\define\do{\dots}
\define\la{\langle}
\define\bsl{\backslash}

\define\lra{\leftrightarrow}
\define\Lra{\Leftrightarrow}

\define\sub{\subset}
\define\bxt{\boxtimes}
\define\T{\times}
\define\ti{\tilde}
\define\nl{\newline}
\redefine\i{^{-1}}

\define\ot{\otimes}
\define\bbq{\bar{\QQ}_l}

\define\Ad{\text{\rm Ad}}

\define\supp{\text{\rm supp}}

\define\a{\alpha}
\redefine\b{\beta}
\redefine\c{\chi}
\define\g{\gamma}
\redefine\d{\delta}
\define\e{\epsilon}

\define\io{\iota}

\define\p{\pi}
\define\ph{\phi}
\define\ps{\psi}
\define\r{\rho}

\redefine\t{\tau}

\define\k{\kappa}
\redefine\l{\lambda}
\define\z{\zeta}
\define\x{\xi}

\define\vt{\vartheta}

\define\Ps{\Psi}

\define\kk{\bold k}

\redefine\tt{\bold t}

\redefine\xx{\bold x}

\define\FF{\bold F}

\define\NN{\bold N}

\define\QQ{\bold Q}

\define\ZZ{\bold Z}

\define\cb{\Cal B}
\define\cc{\Cal C}
\define\cd{\Cal D}
\define\ce{\Cal E}
\define\cf{\Cal F}

\define\cl{\Cal L}

\define\cp{\Cal P}

\define\cs{\Cal S}
\define\ct{\Cal T}

\define\cv{\Cal V}

\define\cz{\Cal Z}

\define\ta{\ti a}

\define\tf{\ti f}

\define\tu{\ti u}

\define\tw{\ti w}

\define\tB{\ti B}

\define\tJ{\ti J}

\define\tP{\ti P}

\define\tX{\ti X}

\define\tZ{\ti Z}

\define\sh{\sharp}

\define\sps{\supset}
\define\BBD{BBD}
\define\BE{B}
\define\DL{DL}
\define\GI{G}
\define\GP{L1}
\define\FU{L2}
\define\CS{L3}
\define\CSS{L4}
\define\CSSS{L5}
\define\IC{L6}
\define\GF{L7}
\define\AIS{L8}

\head Introduction\endhead
\subhead 0.1\endsubhead
Let $G$ be a connected reductive algebraic group over an algebraically closed field 
$\kk$. Let $Z$ be the variety (of the same dimension as $G$) consisting of all pairs 
$(P,gU_P)$ where $P$ runs through a fixed conjugacy class $\cp$ of parabolics of $G$ 
and $gU_P\in G/U_P$ ($U_P$ is the unipotent radical of $P$). Now $G$ acts on $Z$ by 
$h:(P,gU_P)\m(hPh\i,hgh\i U_{hPh\i})$. In this paper we study a class of simple 
$G$-equivariant perverse sheaves on $Z$ which we call "parabolic character sheaves". In
the case where $\cp=\{G\}$, these are the character sheaves on $Z=G$, in the sense of 
\cite{\CS}. In general, they have several common properties with the character sheaves:
they are defined in arbitrary characteristic, and in the case where $\kk$ is an 
algebraic closure of a finite field $\FF_q$ and $G,\cp$ are defined over $\FF_q$, 

(a) {\it the characteristic functions of the parabolic character sheaves that are 
defined over $\FF_q$ form a basis of the space of $G(\FF_q)$-invariant functions on} 
$Z(\FF_q)$
\nl
(under a mild assumption on the characteristic of $\kk$).

\subhead 0.2\endsubhead
We now review the contents of this paper in more detail.

In Section 1 we describe a partition of a partial flag manifold of $G$ (defined over
$\FF_q$) into pieces that are $G(\FF_q)$-stable, generalizing the partition \cite{\DL}
of the full flag manifold. This partition has been considered by the author in 1977
(unpublished); the associated combinatorics has been developed by B\'edard \cite{\BE}
(see Section 2). Although this partition is not needed in our theory of parabolic 
character sheaves, it serves as motivation for it. In Section 3 we define a partition 
of $Z$ into pieces which is governed by the same combinatorics as that in Section 2. In
Section 4 we define the "parabolic character sheaves" on $Z$ in two apparently 
different ways. The first one uses the partition of $Z$ mentioned above and the usual 
notion of character sheaf on a connected component of a possibly disconnected reductive
group. The second one imitates the definition of character sheaves in \cite{\CS}. One 
of our main results is that these two definitions coincide (4.13, 4.16). In Section 4, 
we define a map from the set of parabolic character sheaves on $Z$ to the set of "tame"
local systems on a maximal torus of $G$ modulo the action of a certain subgroup of the 
Weyl group (this extends a known property of character sheaves). In Section 5 we prove 
property 0.1(a). The main result of Section 6 is Corollary 6.8. This has the following 
consequence. Assume that $G$ and $\cp$ (in 0.1) are defined over $\FF_q$. Let $\r$ be 
the character of an irreducible representation of $G(\FF_q)$ over $\bbq$. By summing 
$\r$ over each fibre of $G(\FF_q)@>>>(G/U_P)(\FF_q)$ for $P\in\cp(\FF_q)$, we obtain a 
$G(\FF_q)$-invariant function $Z(\FF_q)@>>>\bbq$. We would like to express this 
function as a linear combination of a small number of elements in a fixed basis of the 
space of all $G(\FF_q)$-invariant functions on $Z(\FF_q)$. (Here "small" means "bounded
by a number independent of $q$".) The basis formed by the characteristic functions of 
the $G(\FF_q)$-orbits on $Z(\FF_q)$ does not have this property. However, the basis 
0.1(a) defined by the parabolic character sheaves does have the required property 
(under a mild assumption) as a consequence of Corollary 6.8. 

\subhead 0.3\endsubhead
The theory of parabolic character sheaves on $Z$ continues to make sense with only
minor modifications when $G$ and $\cp$ (with $G$ simply connected) are replaced by 
$G(\kk((x)))$, $x$ an indeterminate, and a class of parahoric subgroups. (There is no 
known theory of character sheaves on $G(\kk((x)))$ itself.) I believe that this theory 
is a necessary ingredient for establishing the (conjectural) geometric interpretation 
of affine Hecke algebras with unequal parameters proposed in \cite{\AIS}.

\head Contents\endhead
1. A partition of a partial flag manifold.

2. Results of B\'edard.

3. The variety $Z_{J,\d}$ and its partition.

4. Parabolic character sheaves on $Z_{J,\d}$.

5. Central character.

6. The functors $f^J_{J'},e^{J'}_J$.

7. Characteristic functions.

\head 1. A partition of a partial flag manifold\endhead
\subhead 1.1\endsubhead
For any affine connected algebraic group $H$ let $U_H$ be the unipotent radical of $H$.
Let $G$ be as in 0.1. If $P,Q$ are parabolics of $G$ then so is
$$P^Q=(P\cap Q)U_P$$
and we have
$$U_{(P^Q)}=U_P(P\cap U_Q).$$

\subhead 1.2\endsubhead
Assume now that $\kk$ is an algebraic closure of the finite field $\FF_q$ and that we 
are given an $\FF_q$-rational structure on $G$; let $F:G@>>>G$ be the corresponding 
Frobenius map. Let $\cp$ be a conjugacy class of parabolics in $G$. To any $Q\in\cp$ we
attach a sequence of parabolics ${}^0Q\sps{}^1Q\sps{}^2Q\sps\do$ by
$${}^0Q=Q, {}^nQ=({}^{n-1}Q)^{F\i({}^{n-1}Q)} \text{ for } n\ge 1.$$
We say that $Q,Q'\in\cp$ are equivalent if for any $n\ge 0$ there exists an element of
$G$ that conjugates ${}^nQ$ to ${}^nQ'$ and $F({}^nQ)$ to $F({}^nQ')$. The equivalence 
classes form a partition of $\cp$ into finitely many $G^F$-stable subvarieties, or 
"pieces". (In the case where $\cp$ is the variety of Borels, the pieces above are 
precisely the subvarieties $X(w)$ (see \cite{\DL, 1.4}) of the flag manifold.) In the 
general case, using results in \cite{\DL},\cite{\FU} one can show that each piece has 
Euler characteristic in $1+q\ZZ$. Hence

(a) {\it the number of pieces in $\cp$ is equal to the Euler characteristic of $\cp$}
\nl
(at least for $q\gg 0$). Based on this observation, around 1979, I asked R. B\'edard 
(my Ph.D. student at the time) to find a combinatorial explanation for the equality
(a). B\'edard's solution of this problem is reviewed in Section 2.

\subhead 1.3\endsubhead
Let $V$ be a vector space of dimension $n\ge 2$ over $\kk$ (as in 1.2) with a fixed 
$\FF_q$-rational structure with corresponding Frobenius map $F:V@>>>V$. Then $G=GL(V)$ 
is as in 1.2. We identify parabolics in $G$ with partial flags in $V$ in the usual way.
Following \cite{\DL, 2.3} we partition the set of lines in $V$ into locally 
subvarieties $X_1,X_2,\do,X_n$ where $X_j$ is the set of all lines $L$ in $V$ such that
$\dim(L+F(L)+F^2(L)+\do)=j$. The pieces $X_j$ are special cases of the pieces in 1.2.

\subhead 1.4\endsubhead
Let $V$ be a vector space of dimension $2n\ge 4$ over $\kk$ (as in 1.2) with a given
non-degenerate symplectic form $\la,\ra:V\T V@>>>\kk$ and with a fixed $\FF_q$-rational
structure with corresponding Frobenius map $F:V@>>>V$ such that
$\la F(x),F(y)\ra=\la x,y\ra^q$ for all $x,y\in V$. Then the symplectic group $G=Sp(V)$
is as in 1.2. We identify parabolics in $G$ with partial isotropic flags in $V$ in the 
usual way. Following \cite{\GP} we partition the set of lines in $V$ into locally 
closed subvarieties $X_1,X_2,\do,X_n,X'_n,\do,X'_2,X'_1$ where $X_j$ is the set of all 
lines $L$ in $V$ such that $L+F(L)+F^2(L)+\do$ is an isotropic subspace of dimension 
$j$ and $X'_j$ is the set of all lines $L$ in $V$ such that
$$\la L,F(L)\ra=\la L,F^2(L)\ra=\do=\la L,F^{j-1}(L)\ra=0,\la L,F^j(L)\ra\ne 0.$$
The pieces $X_j,X'_j$ are special cases of the pieces in 1.2.

\head 2. Results of B\'edard\endhead
\subhead 2.1\endsubhead
Let $W$ be a Coxeter group, let $I$ be the set of simple reflections and let 
$l:W@>>>\NN$ be the length function. Let $\e:W@>\si>>W$ be a group isomorphism.

In this section we reformulate results of B\'edard \cite{\BE} in a more general 
setting. (In \cite{\BE} it is assumed that $W$ is finite and $\e(I)=I$; we do not 
assume this.)

If $A,A'$ are subsets of a group and $g$ is an element of that group we set 
${}^gA=\Ad(g)A=gAg\i$, $N_A(A')=\{h\in A';{}^hA=A\}$.

For a subset $X$ of $W$ we write $x=\min(X)$ if $x\in X$ and $l(x)<l(x')$ for all 
$x'\in X-\{x\}$. For $J\sub I$ let $W_J$ be the subgroup of $W$ generated by $J$. For 
$J,K\sub I$ let 
$$\align&{}^KW=\{w\in W;w=\min(W_Kw)\},W^J=\{w\in W;w=\min(wW_J)\},\\&
{}^KW^J={}^KW\cap W^J.\endalign$$ 
We recall three known results.

(a) If $x\in{}^KW^{K'},u\in{}^{K'\cap\Ad(x\i)K}W,u\in W_{K'}$, then $xu\in{}^KW$ and 
$l(xu)=l(x)+l(u)$.

(b) If $x\in{}^KW^{K'},x'\in W_KxW_{K'},x'\in{}^KW$, then $x'=xu$ where
$u\in W_{K'}$, $u\in{}^{K'\cap\Ad(x\i)K}W$.

(c) If $x\in W^K,x'\in W_K,K'\sub K$, then $x'\in W^{K'}\Lra xx'\in W^{K'}$.

\subhead 2.2\endsubhead
Let $J,J'\sub I$ be such that $\e(J)=J'$. Let $\ct(J,\e)$ be the set of all sequences 
$(J_n,w_n)_{n\ge 0}$ where $J_n\sub I$ and $w_n\in W$ are such that
$$J=J_0\sps J_1\sps J_2\sps\do,\tag a$$
$$J_n=J_{n-1}\cap\e\i\Ad(w_{n-1})J_{n-1},\text{ for } n\ge 1,\tag b$$
$$w_n\in{}^{\e(J_n)}W^{J_n},\text{ for } n\ge 0,\tag c$$
$$w_n\in W_{\e(J_n)}w_{n-1}W_{J_{n-1}}, \text{ for } n\ge 1.\tag d$$

\subhead 2.3\endsubhead
In the setup of 2.2, let $S(J,\e)$ be the set of all sequences 
$(J_n,J'_n,u_n)_{n\ge 0}$ where $J_n,J'_n$ are subsets of $I$ and $u_n$ are elements of
$W$ such that
$$J=J_0\sps J_1\sps J_2\sps\do,\tag a$$
$$\e(J_n)=\e(J_{n-1})\cap\Ad(u_0u_1\do u_{n-1})J_{n-1}, \text{ for } n\ge 1,\tag b$$
$$J'_0=J',J'_n=J_{n-1}\cap\Ad(u_0u_1\do u_{n-1})\i\e(J_{n-1}), \text{ for } n\ge 1,
\tag c$$
$$u_n\in W_{J_{n-1}}, \text{ for } n\ge 1,\tag d$$
$$u_n\in{}^{J'_n}W^{J_n}, \text{ for } n\ge 0.\tag e$$
From (b),(c) we see that
$$\e(J_n)=\Ad(u_0u_1\do u_{n-1})J'_n, \text{ for } n\ge 1.\tag f$$

\proclaim{Proposition 2.4} There is a unique bijection $S(J,\e)@>\si>>\ct(J,\e)$ such 
that \lb $(J_n,J'_n,u_n)_{n\ge 0}\m(J_n,u_0u_1\do u_n)_{n\ge 0}$.
\endproclaim
Let $(J_n,w_n)_{n\ge 0}\in\ct(J,\e)$. For $n\ge 1$ we set 
$J'_n=J_{n-1}\cap\Ad(w_{n-1}\i)\e(J_{n-1})$. We set $J'_0=J',u_0=w_0$. Then
$u_0\in{}^{J'_0}W^{J_0}$. Let $n\ge 1$. We can find 
$v\in W_{\e(J_{n-1})}$ such that $v\i w_n\in{}^{\e(J_{n-1})}W$. We have
$$v\i w_n\in W_{\e(J_{n-1})}w_{n-1}W_{J_{n-1}},w_{n-1}\in{}^{\e(J_{n-1})}W^{J_{n-1}}.$$
Hence, by 2.1(b) and 2.1(a) we have $v\i w_n=w_{n-1}u_n$ with $u_n\in W_{J_{n-1}}$,
$u_n\in{}^{J'_n}W$, $l(w_{n-1}u_n)=l(w_{n-1})+l(u_n)$. From 2.2(d) we deduce that 
$vw_{n-1}\in W_{\e(J_n)}w_{n-1}W_{J_{n-1}}$. From 2.2(b) we have 
$W_{\e(J_n)}w_{n-1}\sub w_{n-1}W_{J_{n-1}}$ hence $vw_{n-1}\in w_{n-1}W_{J_{n-1}}$ and
$v\in w_{n-1}W_{J_{n-1}}w_{n-1}\i\cap W_{\e(J_{n-1})}=W_{\e(J_n)}$. From 
$v\in W_{\e(J_n)},w_n\in{}^{\e(J_n)}W,v\i w_n\in{}^{\e(J_n)}W$, we deduce that $v=1$.
Hence $w_n=w_{n-1}u_n$. Let $v'\in W_{J_n}$. Since $w_n\in W^{J_n}$, we have 
$$l(w_nv')=l(w_n)+l(v')=l(w_{n-1}u_n)+l(v')=l(w_{n-1})+l(u_n)+l(v').$$
Since $w_{n-1}\in W^{J_{n-1}}$ and $u_nv'\in W_{J_{n-1}}$, we have 
$l(w_nv')=l(w_{n-1}u_nv')=l(w_{n-1})+l(u_nv')$. Thus,
$l(w_{n-1})+l(u_n)+l(v')=l(w_{n-1})+l(u_nv')$ hence $l(u_n)+l(v')=l(u_nv')$. We see 
that $u_n\in W^{J_n}$. Thus, $u_n\in{}^{J'_n}W^{J_n}$ and 
$(J'_n,J_n,u_n)_{n\ge 0}\in S(J,\e)$. Thus, we have a map 
$\ct(J,\e)@>>>S(J,\e)$, $(J_n,w_n)_{n\ge 0}\m(J_n,J'_n,u_n)_{n\ge 0}$ where 
$u_0=w_0,u_n=w_{n-1}\i w_n$ for $n\ge 1$. We now construct an inverse of this map. Let 
$(J_n,J'_n,u_n)_{n\ge 0}\in S(J,\e)$. 

(a) {\it For any $0\le k\le n$ we have $u_ku_{k+1}\do u_n\in{}^{J'_k}W$ and} 
$l(u_ku_{k+1}\do u_n)=l(u_k)+l(u_{k+1})+\do+l(u_n)$.
\nl
We use descending induction on $k\in[0,n]$. For $k=n$, (a) is clear. Assume now that 
$k<n$ and that (a) is true when $k$ is replaced by $k+1$. We have
$$u_k\in{}^{J'_k}W^{J_k},u_{k+1}\do u_n\in{}^{J'_{k+1}}W={}^{J_k\cap\Ad(u_k\i)J'_k}W,$$
$u_{k+1}\do u_n\in W_{J_k}$ and $l(u_{k+1}\do u_n)=l(u_{k+1})+\do+l(u_n)$. Using 2.1(a)
we deduce $u_ku_{k+1}\do u_n\in{}^{J'_k}W$ and 
$$l(u_ku_{k+1}\do u_n)=l(u_k)+l(u_{k+1}\do u_n)=l(u_k)+l(u_{k+1})+\do+l(u_n).$$
(a) is proved.

(b) {\it For any $n\ge 0$ we have $u_0u_1\do u_n\in{}^{J'}W$.}
\nl
This is a special case of (a).

(c) {\it For any $n\ge 0$ we have $u_0u_1\do u_n\in W^{J_n}$.}
\nl
We use induction on $n$. For $n=0$, (c) is clear. Assume now that $n>0$ and that (c) is
true when $n$ is replaced by $n-1$. We have 
$$u_0u_1\do u_{n-1}\in W^{J_{n-1}},u_n\in W^{J_n},u_n\in W_{J_{n-1}},J_n\sub J_{n-1}.$$
Using 2.1(c) we deduce $u_0u_1\do u_{n-1}u_n\in W^{J_n}$. This proves (c).

(d) {\it For any $n\ge 0$ we have $u_0u_1\do u_n\in{}^{J'}W^{J_n}$.}
\nl
This follows from (b) and (c).

(e) {\it For any $n\ge 0$ we have $u_0u_1\do u_n\in{}^{\e(J_n)}W^{J_n}$.}
\nl
This follows from (d) since $\e(J_n)\sub J'_0$. 

Using (e), define $S(J,\e)@>>>\ct(J,\e)$ by 
$(J_n,J'_n,u_n)_{n\ge 0}\m(J_n,u_0u_1\do u_n)_{n\ge 0}$. Clearly, we have defined two 
inverse bijections between $S(J,\e),\ct(J,\e)$. The proposition is proved.

\proclaim{Proposition 2.5} $(J_n,J'_n,u_n)_{n\ge 0}\m u_0u_1\do u_m$ for $m\gg 0$ is a
well defined bijection $\ph:S(J,\e)@>>>{}^{J'}W$.
\endproclaim
Let $(J_n,w_n)_{n\ge 0}\in\ct(J,\e)$. There exists $n_0\ge 1$ such that $J_{n-1}=J_n$ 
for $n\ge n_0$. For such $n$ we have 
$$w_n\in{}^{\e(J_n)}W^{J_n},w_{n-1}\in{}^{\e(J_n)}W^{J_n},
w_n\in W_{\e(J_n)}w_{n-1}W_{J_n}$$
hence $w_n=w_{n-1}$. Thus, if $(J_n,J'_n,u_n)_{n\ge 0}\in S(J,\e)$ then for $n\gg 0$ we
have $u_n=1$. Hence there is a well defined element $w\in W$ such that 
$w=u_0u_1\do u_n$ for $n\gg 0$. By 2.4(d), we have $w\in{}^{J'}W$. Hence $\ph$ is a
well defined map.

(a) {\it For any $n\ge 0$ we have $u_0u_1\do u_n=\min(W_{J'}wW_{J_n})$.}
\nl
By 2.4(d), we have $u_0u_1\do u_n\in{}^{J'}W^{J_n}$. Hence it suffices to show that 
$w\in u_0u_1\do u_nW_{J_n}$. Take $N>n$ such that $u_0u_1\do u_N=w$. Then 
$$w=(u_0\do u_n)u_{n+1}u_{n+2}\do u_N$$ 
and it suffices to show that $u_{n+1}u_{n+2}\do u_N\in W_{J_n}$. This follows from 
$$u_{n+1}\in W_{J_n},u_{n+2}\in W_{J_{n+1}}\sub W_{J_n},\do,u_N\in W_{J_{N-1}}\sub 
W_{J_n}.$$ 
This proves (a).

We show that $\ph$ is a bijection. Assume that the images $w,\tw$ of 
$$(J_n,J'_n,u_n)_{n\ge 0}\in S(J,\e),\qua(\tJ_n,\tJ'_n,\tu_n)_{n\ge 0}\in S(J,\e)$$
under $\ph$ satisfy $w=\tw$. We show by induction on $n\ge 0$ that

(b) $J'_k=\tJ'_k,J_k=\tJ_k,u_k=\tu_k$ for $k\in[0,n]$.
\nl
For $n=0$ this holds since $J_0=\tJ_0=J,J'_0=\tJ'_0=J'$ and 
$u_0=\tu_0=\min(W_{J'_0}wW_{J_0})$ (see (a)). Assume now that $n>0$ and that (b) holds
when $n$ is replaced by $n-1$. From 2.3(b),(c) we deduce that $J_n=\tJ_n,J'_n=\tJ'_n$. 
From (a) we have
$$u_0u_1\do u_n=\min(W_{J'_0}wW_{J_n})=\min(W_{J'_0}wW_{\tJ_n})=\tu_0\tu_1\do\tu_n
=u_0u_1\do u_{n-1}\tu_n$$
hence $u_n=\tu_n$. Thus (b) holds. We see that 
$(J_n,J'_n,u_n)_{n\ge 0}=(\tJ_n,\tJ'_n,\tu_n)_{n\ge 0}$. Thus, $\ph$ is injective. 

We define an inverse to $\ph$. Let $w\in{}^{J'}W$. We define by induction on $n\ge 0$ 
a sequence $(J_n,J'_n,u_n)_{n\ge 0}$ as follows. We set $J_0=J,J'_0=J'$,
$u_0=\min(W_{J'_0}wW_{J_0})$. Assume that $n>0$ and that $J'_k,J_k,u_k$ are defined for
$k\in[0,n-1]$. We define subsets $J_n,J'_n$ of $J_{n-1}$ by
$$\align&\e(J_n)=\e(J_{n-1})\cap\Ad(u_0u_1\do u_{n-1})J_{n-1},\\&
J'_n=J_{n-1}\cap\Ad(u_0u_1\do u_{n-1})\i\e(J_{n-1})\endalign$$
and we define $u_n$ by $u_0u_1\do u_n=\min(W_{J'}wW_{J_n})$. This completes the
inductive definition. Using $\e(J_{n-1})\sub J'$, we see that
$$J'_n\sub\Ad(u_0u_1\do u_{n-1})\i J' \text{ for } n\ge 1.\tag c$$
We show that
$$u_n\in W_{J_{n-1}} \text{ for } n\ge 1.\tag d$$
Now $u_0u_1\do u_n\in W_{J'}wW_{J_n}\sub W_{J'}wW_{J_{n-1}}$ and
$u_0u_1\do u_{n-1}=\min(W_{J'}wW_{J_{n-1}})$. Moreover, $u_0u_1\do u_n\in{}^{J'}W$. By
2.1(b) we have $u_0u_1\do u_n=(u_0u_1\do u_{n-1})u$ with 
$$u\in W_{J_{n-1}},u\in{}^{J_{n-1}\cap\Ad(u_0u_1\do u_{n-1})\i J'}W.$$
Thus $u=u_n$ and (d) follows. We show by induction on $n\ge 0$ that
$$u_n\in{}^{J'_n}W^{J_n}.\tag e$$
For $n=0$ this is clear. Assume now that $n>0$ and that (e) holds when $n$ is replaced 
by $n-1$. By the argument in the proof of (d) we have 
$u_n\in{}^{J_{n-1}\cap\Ad(u_0u_1\do u_{n-1})\i J'}W$. Since $J'_n\sub J_{n-1}$, we see
from (c) that 
$$J'_n\sub J_{n-1}\cap\Ad(u_0u_1\do u_{n-1})\i J'.$$
Hence $u_n\in{}^{J'_n}W$. Now 
$$u_0u_1\do u_{n-1}\in W^{J_{n-1}},u_n\in W_{J_{n-1}},
(u_0u_1\do u_{n-1})u_n\in W^{J_n},J_n\sub J_{n-1}.$$
Using 2.1(c) we deduce that $u_n\in W^{J_n}$. Combining this with $u_n\in{}^{J'_n}W$
gives $u_n\in{}^{J'_n}W^{J_n}$. Thus (e) is established. We see that 
$(J_n,J'_n,u_n)_{n\ge 0}\in S(J,\e)$. We show that

(f) {\it if $n\gg 0$, then } $u_0u_1\do u_n=w.$
\nl
For any $n\ge 0$ we have $u_0u_1\do u_n=\min(W_{J'}wW_{J_n})$. Also, $w\in{}^{J'}W$.
Using 2.1(b) we have $w=(u_0u_1\do u_n)u$ with 
$$u\in W_{J_n},u\in{}^{J_n\cap\Ad(u_0u_1\do u_n)\i J'}W.$$ 
By (c), we have $J'_{n+1}\sub J_n\cap\Ad(u_0u_1\do u_n)\i J'$. Hence 
$u\in{}^{J'_{n+1}}W$. Now $J'_{n+1}\sub J_n$. Assuming that $n\gg 0$, we have 
$J_n=J_{n+1}$. From 2.3(f) we see that $\sh J'_{n+1}=\sh J_{n+1}$, hence 
$\sh J'_{n+1}=\sh J_n$. Hence the inclusion $J'_{n+1}\sub J_n$ must be an equality. We 
see that $u\in{}^{J_n}W$. This, combined with $u\in W_{J_n}$ implies $u=1$ and proves 
(f). Thus, we have defined $\ps:{}^{J'}W@>>>S(J,\e)$ such that $\ph\ps=1$. Hence $\ph$ 
is bijective. 

\subhead 2.6\endsubhead
Let $(J_n,J'_n,u_n)_{n\ge 0}\in S(J,\e)$ and let $w$ be its image under $\ph$. For 
$n\gg 0$ we have $u_0u_1\do u_{n-1}=w$ hence $\Ad(w)J'_n=\e(J_n)$. By the proof of
2.5(f), for $n\gg 0$ we have also $J_n=J'_n=J_{n+1}=J'_{n+1}=\do$. Thus there is a well
defined subset $J_\iy$ of $J$ such that $J_n=J'_n=J_\iy$ for $n\gg 0$ and 
$\Ad(w)J_\iy=\e(J_\iy)$.

\subhead 2.7\endsubhead
Let $G$ be as in 0.1. Let $\cb$ be the variety of Borel subgroups of $G$. Let $W$ be 
the set of $G$-orbits on $\cb\T\cb$ ($G$ acts by conjugation on both factors). Then $W$
is naturally a finite Coxeter group; let $I$ be the set of simple reflections (the 
$G$-orbits of dimension $\dim\cb+1$). For $B,B'\in\cb,w\in W$ we write $\po(B,B')=w$ if
the $G$-orbit of $(B,B')$ is $w$. If $P$ is a parabolic of $G$, the set of all $w\in W$
such that $w=\po(B,B')$ for some $B,B'\in\cb,B\sub P,B'\sub P$, is of the form $W_J$ 
(as in 2.1) for a well defined subset $J\sub I$. We then say that $P$ has type $J$. For
$J\sub I$, let $\cp_J$ be the set of all parabolics of type $J$ of $G$.

For $P\in\cp_J,Q\in\cp_K$ there is a well defined element $u=\po(P,Q)\in{}^JW^K$ such 
that $\po(B,B')\ge u$ (standard partial order on $W$) for any 
$B,B'\in\cb,B\sub P,B'\sub Q$ and $\po(B_1,B'_1)=u$ for some 
$B_1,B'_1\in\cb,B_1\sub P,B'_1\sub Q$; we then have $B_1\sub P^Q,B'_1\sub Q^P$. We have
$$P^Q\in\cp_{J\cap\Ad(u)K}.\tag a$$
Now $(P,Q)\m u$ defines a bijection between the set of $G$-orbits on $\cp_J\T\cp_K$ and
${}^JW^K$.

\subhead 2.8\endsubhead
Let $G,F$ be as in 1.2. Let $W,I$ be as in 2.7. The bijection $\d:W@>>>W$ induced by 
$F:G@>>>G$ satisfies $\d(I)=I$. Let $J\sub I$. We show that the pieces of $\cp_J$ 
defined in 1.2 are naturally indexed by $\ct(J,\d)$.

Let $Q\in\cp_{J}$. To $Q$ we associate a sequence $(J_n,w_n)_{n\ge 0}$ with 
$J_n\sub I$, $w_n\in W$ and a sequence $({}^nQ)_{n\ge 0}$ with ${}^nQ\in\cp_{J_n}$. We
set 
$${}^0Q=Q,J_0=J,w_0=\po(F({}^0Q),{}^0Q).$$
Assume that $n\ge 1$, that ${}^mQ,J_m,w_m$ are already defined for $m<n$ and that 
$w_m=\po(F({}^mQ),{}^mQ),{}^mQ\in\cp_{J_m}$ for $m<n$. Let
$$J_n=J_{n-1}\cap\d\i\Ad(w_{n-1})J_{n-1},$$
$${}^nQ=({}^{n-1}Q)^{F\i({}^{n-1}Q)}\in\cp_{J_n},$$
$$w_n=\po(F({}^nQ),{}^nQ)\in{}^{\d(J_n)}W^{J_n}.$$
This completes the inductive definition. From Lemma 3.2(c) (with $P,P',Z$ replaced by 
${}^{n-1}Q,F({}^{n-1}Q),{}^nQ$) we see that $w_n\in w_{n-1}W_{J_{n-1}}$ for $n\ge 1$. 
Thus, $(J_n,w_n)_{n\ge 0}\in\ct(J,\d)$. Thus, $Q\m(J_n,w_n)_{n\ge 0}$ is a map
$\cp_J@>>>\ct(J,\d)$. The fibre of this map at $\tt\in\ct(J,\d)$ is denoted by 
${}^\tt\cp_J$. Clearly, $({}^\tt\cp_J)_{\tt\in\ct(J,\d)}$ is a partition of $\cp_J$ 
into locally closed subvarieties (the same as in 1.2). The $G^F$-action on $\cp_J$ 
given by $h:Q\m{}^hQ$ preserves each of the pieces ${}^\tt\cp_J$. 

\head 3. The variety $Z_{J,\d}$ and its partition\endhead
\subhead 3.1\endsubhead
Let $G$ be as in 0.1. Let $W,I$ be as in 2.7. Let $\hG$ be a possibly disconnected 
reductive algebraic group over $\kk$ with identity component $G$ and let $G^1$ be a 
fixed connected component of $\hG$. There is a unique isomorphism $\d:W@>\si>>W$ such 
that $\d(I)=I$ and
$$P\in\cp_J,g\in G^1\imp{}^gP\in\cp_{\d(J)}.$$
Let $P\in\cp_J,Q\in\cp_K,u=\po(P,Q)$. We say that $P,Q$ are {\it in good position} if
they have a common Levi or, equivalently, if $\Ad(u\i)(J)=K$. In this case we have 
$P^Q=P,Q^P=Q$.

We fix $J\sub I$. 

\proclaim{Lemma 3.2} Let $P\in\cp_J,P'\in\cp_{J'}$, $J'\sub I$. Let $a=\po(P',P)$. Let 
$X=P'{}^P,Y=P^{P'}$ and let $Z$ be a parabolic subgroup of $P$. Let $b=\po(Y,Z)$.

(a) Let $Y'$ be a parabolic subgroup of $P$ of the same type as $Y$ such that $X,Y'$ 
are in good position and $\po(X,Y')=a$. Then $Y'=Y=P^X$.

(b) $X$ contains a Levi of $Y\cap Z$.

(c) $\po(X,Z)=ab$.

(d) $X^{(Y^Z)}=X^Z$.
\endproclaim
We prove (a). Since $X,Y'$ are in good position, for any Borel $B'$ in $Y'$ there 
exists a Borel $B$ in $X$ such that $\po(B,B')=a$. Since $\po(X,P)=a$ and $B'\sub P$ we
have $B'\sub P^X$. Since $Y'$ is the union of its Borels, we have $Y'\sub P^X$. Now 
$Y'$ is of type $J\cap\Ad(a\i)J'$, and $P^X$ is of type 
$J\cap\Ad(a\i)(J'\cap\Ad(a)J)=J\cap\Ad(a\i)J'$. Thus $Y',P^X$ have the same type, hence
$Y'=P^X$. Replacing in this argument $Y'$ by $Y$ we obtain $Y=P^X$. 

We prove (b). Since $Z\sub P$ we have $U_P\sub Z$. Hence
$Y\cap Z=P^{P'}\cap Z=((P\cap P')U_P)\cap Z=(P\cap P'\cap Z)U_P$. Now 
$U_P\sub U_{Y\cap Z}$. Hence if $L$ is a Levi of $P\cap P'\cap Z$ then $L$ is a Levi of
$Y\cap Q$. Now $L\sub P\cap P'\sub X$.

We prove (c). Let $\ta=\po(X,Z)$. We have $X\in\cp_K$ with $K=J'\cap\Ad(a)J$. We have 
$a\in{}^{J'}W^J$ and $K\sub J'$ hence $a\in{}^KW^J$. From the definition of $\ta$ we 
have $\ta\in{}^KW$. Since $Z\sub P$, we have $\ta\in W_KaW_J$. Using 2.1(c) with $x,x'$
replaced by $a,\ta$, we see that $\ta=av$ with $v\in W_J$. Let $B,B'\in\cb$ be such 
that $B\sub X$, $B'\sub Z$, $\po(B,B')=\ta$. Since $a\in W^J$, we have 
$l(av)=l(a)+l(v)$ and there is a unique $B''\in\cb$ such that $\po(B,B'')=a$, 
$\po(B'',B')=v$. Since $B'\sub P$ and $\po(B'',B)\in W_J$, we have $B''\sub P$. Since
$B\sub P',B''\sub P$ and $\po(B,B'')=\po(P',P)=a$, we have $B''\sub P^{P'}=Y$. Since 
$B''\sub Y,B'\sub Z$, we have $v\ge b$. We can find $B_1,B_2\in\cb$ such that 
$B_1\sub Y$, $B_2\sub Z$, $\po(B_1,B_2)=b$. Since $\po(P',P)=a$ and $B_1\sub P^{P'}$, 
we can find $B_0\in\cb$ such that $B_0\sub P'$, $\po(B_0,B_1)=a$. Since 
$a\in W^J,b\in W_J$, we have $\po(B_0,B_2)=ab$. We have $B_0\sub P'{}^P=X$, $B_2\sub Z$
hence $\po(B_0,B_2)\ge\po(X,Z)$, that is, $ab\ge\ta=av$. Thus, we have $v\ge b$ and 
$ab\ge av$. Since $a\in W^J$ and $b,v\in W_J$ we have $b=v$. Thus, $\ta=ab$. 

We prove (d). Let $B$ be a Borel in $X^Z$. Using (c), we can find a Borel $B'$ in $Z$
such that $\po(B,B')=ab$. Let $B''$ be the unique Borel such that $\po(B,B'')=a$,
$\po(B'',B')=b$. Now $B''\sub P$ since $B'\sub P$ and $\po(B'',B')=b\in W_J$. Since 
$\po(X,P)=\po(B,B'')=a$, we have $B''\sub P^X$ hence $B''\sub Y$ (see (a)). Since 
$\po(Y,Z)=\po(B'',B')=b$, we have $B''\sub Y^Z$. Since $B\sub X,B''\sub Y^Z$ and 
$\po(B,B'')=a=\po(P',P)$ where $X\sub P',Y^Z\sub P$, we have $\po(B,B'')=a=\po(X,Y^Z)$.
Hence $B\sub X^{(Y^Z)}$. Thus any Borel in $X^Z$ is contained in $X^{(Y^Z)}$. Since 
$X^Z$ is the union of its Borels, we have $X^Z\sub X^{(Y^Z)}$. To show that 
$X^Z=X^{(Y^Z)}$, it suffices to show that $X^Z,X^{(Y^Z)}$ have the same type, or that 
$K\cap\Ad(ab)K''=K\cap\Ad(a)(K'\cap\Ad(b)K'')$ (where $X,Y,Z$ are of type $K,K',K''$) 
or that $\Ad(a)K'=K$ which follows from $K=J'\cap\Ad(a)J$, $K'=J\cap\Ad(a\i)J'$. The 
lemma is proved.

\subhead 3.3\endsubhead
For $(P,P')\in\cp_J\T\cp_{\d(J)}$,
$$A(P,P')=\{g\in G^1;{}^gP=P'\}$$
is a single $P$-orbit (resp. $P'$-orbit) for right (resp. left) translation on $G$. Let
$\cz_{J,\d}$ be the set of all triples $(P,P',g)$ where $P\in\cp_J,P'\in\cp_{\d(J)}$, 
$g\in A(P,P')$. 

Let $(P,P',g)\in\cz_{J,\d}$; let $z=\po(P',P)$. Let $J_1=J\cap\d\i\Ad(z)(J)$. We set
$$P^1=g\i P'{}^Pg\in\cp_{J_1},P'{}^1=P'{}^P\in\cp_{\d(J_1)}.$$
We write $(P^1,P'{}^1)=\a'(P,P',g)$. We have $(P^1,P'{}^1,g)\in\cz_{J_1,\d}$. We have
$P^1\sub P,P'{}^1\sub P'$ and, by Lemma 3.2(c) (with $Z=P^1$), we have
$\po(P'{}^1,P^1)\in zW_J$. 

\proclaim{Lemma 3.4} Let $g,g'\in A(P,P'),u'\in U_{P'},u\in U_P$ be such that 
$g'=u'gu$. Then 

(a) $g'=u'_1gu_1$ where $u'_1\in U_{P'{}^1},u_1\in U_{P^1}$;

(b) we have $\a'(P,P',g')=(P^1,P'{}^1)$.
\endproclaim
To prove (a), we may assume that $u'=1$ (since $U_{P'}g=gU_P$). Since $P^1\sub P$, we 
have $U_P\sub U_{P^1}$ hence $g'\in gU_{P^1}$. This proves (a). To prove (b), we may 
assume that $u=1$ (since $U_{P'}g=gU_P$). Since $U_{P'}\sub P'{}^P$, we have 
$u'\in P'{}^P$ hence ${}^{g'{}\i}(P'{}^P)={}^{g\i}({}^{u'{}\i}(P'{}^P))=P^1$. The lemma
is proved.

\proclaim{Lemma 3.5} Let $g,g'\in A(P,P')$. Assume that 
$\a'(P,P',g)=\a'(P,P',g')=(P^1,P'{}^1)$ and that $g'\in U_{P'{}^1}g=gU_{P^1}$. There
exist $x\in U_P\cap P',w'\in U_{P'}$ such that $g'=w'xg$.
\endproclaim
We have $g'=u'g$ where $u'\in U_{P'{}^P}$. We have $u'=w'x$ where $w'\in U_{P'}$,
$x\in U_P\cap P'$. Hence $g'=w'xg$. 

\proclaim{Lemma 3.6} Let $z\in{}^{\d(J)}W^J$, $J_1=J\cap\d\i\Ad(z)(J)$. Let 

$\cz_{J,\d}^*=\{(P,P',g)\in\cz_{J,\d};\po(P',P)=z\}$,

$\cz_{J_1,\d}^\dag=\{(Q,Q',g)\in\cz_{J_1,\d};\po(Q,Q')\in zW_J\}$.
\nl
Define $f:\cz_{J,\d}^*@>>>\cz_{J_1,\d}^\dag$ by $(P,P',g)\m(P^1,P'{}^1,g)$ where 
$(P^1,P'{}^1)=\a'(P,P',g)$. Then $f$ is an isomorphism.
\endproclaim
We show only that $f$ is bijective. For $(P,P',g)\m\cz_{J,\d}^*$, $P$ (resp. $P'$) is 
the unique parabolic of type $J$ (resp. $\d(J)$) that contains $P^1$ (resp. $P'{}^1$). 
Hence $f$ is injective. We show that $f$ is surjective. Let 
$(Q,Q',g)\in\cz_{J_1,\d}^\dag$. Let $P$ (resp. $P'$) be the unique parabolic of type 
$J$ (resp. $\d(J)$) that contains $Q$ (resp. $Q'$). Clearly, $\po(P',P)=z$ and 
${}^gP=P'$. It suffices to show that $P'{}^P=Q'$. We have $\po(Q',Q)=zu$ where 
$u\in W_J$. We can find $B,B'\in\cb$ such that $B\sub Q,B'\sub Q',\po(B',B)=zu$. Since
$l(zu)=l(z)+l(u)$, we can find $B''\in\cb$ such that $\po(B',B'')=z$, $\po(B'',B)=u$.
Since $u\in W_J$ and $B\sub P$, we have $B''\sub P$. Since $B'\sub P',B''\sub P$, 
$\po(B',B'')=z$, we have $B'\sub P'{}^P$. Since $Q',P'{}^P$ are in $\cp_{\d(J_1)}$ and 
both contain $B'$, we have $Q'=P'{}^P$. The lemma is proved.

\subhead 3.7\endsubhead
We fix $z\in{}^{\d(J)}W^J$. Let $J_1=J\cap\d\i\Ad(z)J$. Let 
$(Q,Q',\g_1)\in Z_{J_1,\d}$ be such that $\po(Q',Q)\in zW_J$. Let $\cf$ be the set of 
all $(P,P',\g)\in Z_{J,\d}$ such that $\po(P',P)=z$, $\g\sub\g_1$ and 
$\a'(P,P',g)=(Q,Q')$ for some/any $g\in\g$ (see 3.4). (Note that $P,P'$ are uniquely 
determined by $Q,Q'$.) By Lemma 3.6 we have $\cf\ne\em$. 

\proclaim{Lemma 3.8} Let $(P,P',\g)\in\cf$. Let $g\in\g$. Then $v\m(P,P',U_{P'}vgU_P)$ 
is a well defined, surjective map $\k:U_P\cap P'@>>>\cf$, 
\endproclaim
Let $(P,P',\g')\in\cf$. Let $g'\in\g'$. By Lemma 3.5 there exist $v\in U_P\cap P'$,
$w'\in U_{P'}$ such that $g'=w'vg$. Then $\g'=U_{P'}vgU_P$. The lemma is proved.

\proclaim{Lemma 3.9} In the setup of Lemma 3.8, the following two conditions for $v,v'$
in $U_P\cap P'$ are equivalent:

(i) $\k(v)=\k(v')$;

(ii) $v'=dv$ for some $d\in U_P\cap U_{P'}$.
\endproclaim
Assume that (i) holds. We have $vg\in U_{P'}v'gU_P=U_{P'}v'U_{P'}g$ hence
$v\in U_{P'}v'U_{P'}=U_{P'}v'$. Thus, $v'=dv$ where $d\in U_{P'}$ and we have
automatically $d\in U_P\cap U_{P'}$. Thus, (ii) holds. The converse is immediate.

\subhead 3.10\endsubhead
We see that $\k$ defines a bijection $(U_P\cap U_{P'})\bsl(U_P\cap P')@>\si>>\cf$. One 
can check that this is an isomorphism of algebraic varieties. Since $U_P\cap P'$ is a 
connected unipotent group and $U_P\cap U_{P'}$ is a connected closed subgroup of it, we
see that 

(a) $\cf$ is isomorphic to an affine space of dimension 
$\dim((U_P\cap U_{P'})\bsl(U_P\cap P'))$.

\subhead 3.11\endsubhead
To any $(P,P',g)\in\cz_{J,\d}$ we associate a sequence $(J_n,w_n)_{n\ge 0}$ with 
$J_n\sub I$, $w_n\in W$ and a sequence $(P^n,P'{}^n,g)_{n\ge 0}$ with 
$(P^n,P'{}^n,g)\in\cz_{J_n,\d}$. We set 
$$P^0=P,P'{}^0=P',J_0=J,w_0=\po(P'{}^0,P^0).$$  
Assume that $n\ge 1$, that $P^m,P'{}^m,J_m,w_m$ are already defined for $m<n$ and that 
$w_m=\po(P'{}^m,P^m),P^m\in\cp_{J_m},P'{}^m={}^gP^m$ for $m<n$. Let
$$J_n=J_{n-1}\cap\d\i\Ad(w_{n-1})(J_{n-1}),$$
$$P^n=g\i(P'{}^{n-1})^{P^{n-1}}g\in\cp_{J_n},\qua
P'{}^n=(P'{}^{n-1})^{P^{n-1}}\in\cp_{\d(J_n)},$$
$$w_n=\po(P'{}^n,P^n)\in{}^{\d(J_n)}W^{J_n}.$$
This completes the inductive definition. From Lemma 3.2(c) (with $P,P',Z$ replaced by 
$P^{n-1},P'{}^{n-1},P^n$) we see that $w_n\in w_{n-1}W_{J_{n-1}}$ for $n\ge 1$. Thus, 
$(J_n,w_n)_{n\ge 0}\in\ct(J,\d)$. We write $(J_n,w_n)_{n\ge 0}=\b'(P,P',g)$. For 
$\tt\in\ct(J,\d)$ let 
$${}^\tt\cz_{J,\d}=\{(P,P',g)\in\cz_{J,\d};\b'(P,P',g)=\tt\},$$
$${}^\tt Z_{J,\d}=\{(P,P',\g)\in Z_{J,d};\b'(P,P',g)=\tt\text{ for some/any } g\in\g\}.
$$
(The equivalence of "some/any" follows from Lemma 3.4.) Clearly, 
$({}^\tt\cz_{J,\d})_{\tt\in\ct(J,\d)}$ is a partition of $\cz_{J,\d}$ into locally 
closed subvarieties and $({}^\tt Z_{J,\d})_{\tt\in\ct(J,\d)}$ is a partition of
$Z_{J,\d}$ into locally closed subvarieties. The $G$-action on $\cz_{J,\d}$ given by
$h:(P,P',g)\m({}^hP,{}^hP',{}^h\g)$ preserves each of the pieces ${}^\tt\cz_{J,\d}$. 
Similarly, the natural action of $G$ on $Z_{J,\d}$ preserves each of the pieces 
${}^\tt Z_{J,\d}$. Clearly, $(P,P',g)\m(P^1,P'{}^1,g)$ is a morphism
$\vt':{}^\tt\cz_{J,\d}@>>>{}^{\tt_1}\cz_{J_1,\d}$ where for 
$\tt=(J_n,w_n)_{n\ge 0}\in\ct(J,\d)$ we set $\tt^1=(J_n,w_n)_{n\ge 1}\in\ct(J_1,\d)$; 
it induces a morphism 
$$\vt:{}^\tt Z_{J,\d}@>>>{}^{\tt_1}Z_{J_1,\d}.$$ 

\proclaim{Lemma 3.12}(a) The morphism $\vt':{}^\tt\cz_{J,\d}@>>>{}^{\tt_1}\cz_{J_1,\d}$
is an isomorphism.

(b) The morphism $\vt:{}^\tt Z_{J,\d}@>>>{}^{\tt_1}Z_{J_1,\d}$ is an affine space
bundle with fibres of dimension $\dim((U_P\cap U_{P'})\bsl(U_P\cap P'))$ where 
$(P,P',\g)$ is any triple in ${}^\tt Z_{J,\d}$.

(c) Consider the map $\bvt$ from the set of $G$-orbits on ${}^\tt Z_{J,\d}$ to the set
of $G$-orbits on ${}^{\tt_1}Z_{J_1,\d}$ induced by $\vt$. Then $\bvt$ is a bijection.
\endproclaim
(a) follows from Lemma 3.6. We prove (b). Now 
$\dim((U_P\cap U_{P'})\bsl(U_P\cap P'))$ in (b) is independent of the choice of 
$(P,P',\g)$; it depends only on $\po(P,P')$ which is constant on ${}^\tt Z_{J,\d}$. 
From 3.10(a) we see that each fibre of $\vt$ is an affine space of the indicated 
dimension. The verification of local triviality is omitted.

We prove (c). Since $\vt$ is surjective (see (b)) and $G$-equivariant, $\bvt$ is well 
defined and surjective. We show that $\bvt$ is injective. Let
$(P,P',\g),(\tP,\tP',\ti\g)$ be two triples in ${}^\tt Z_{J,\d}$ whose images under
$\vt$ are in the same $G$-orbit; we must show that these two triples are in the same 
$G$-orbit. Since $\vt$ is $G$-equivariant, we may assume that 
$\vt(P,P',\g)=\vt(\tP,\tP',\ti\g)=(Q,Q',\g_1)\in{}^{\tt_1}Z_{J_1,\d}$. Define $\cf$ in 
terms of $(Q,Q',\g_1)$ as in 3.7. Then $(P,P',\g)\in\cf,(\tP,\tP',\ti\g)\in\cf$. Since 
$P,\tP$ are parabolics of the same type containing $Q$ we have $P=\tP$. Since $P',\tP'$
are parabolics of the same type containing $Q'$ we have $P'=\tP'$. Let $g\in\g$. By 
Lemma 3.8, we have $\ti\g=U_{P'}vgU_P$ for some $v\in U_P\cap P'$. We have also 
$\tP={}^vP$ (since $v\in P$), $\tP'={}^vP'$ (since $v\in P'$), 
$$\ti\g=U_{P'}vgU_P=vU_{P'}gU_P=vU_{P'}gU_Pv\i=v\g v\i$$
(since $v$ normalizes $U_{P'}$ and $v\in U_P$). Thus, $(\tP,\tP',\ti\g)$ is obtained by
the action of $v\in G$ on $(P,P',\g)$, hence $(\tP,\tP',\ti\g)$ is in the $G$-orbit of 
$(P,P',\g)$. The lemma is proved.

\proclaim{Lemma 3.13} Let $\tt=(J_n,w_n)_{n\ge 0}\in\ct(J,\d)$. Then ${}^\tt Z_{J,\d}$ 
is an iterated affine space bundle over a fibre bundle over $\cp_{J_n}$ with fibres 
isomorphic to $P/U_P$ where $P\in\cp_{J_n}$ ($n\gg 0$). In particular, 
${}^\tt Z_{J,\d}\ne\em$.
\endproclaim
Assume first that $\tt$ is such that $J_n=J$ for all $n$ and $w_n=w$ for all $n$ (here 
$w\in W$). In this case, ${}^\tt Z_{J,\d}$ is the set of all $(P,P',\g)$ where 
$P\in\cp_J$, $P'\in\cp_{\d(J)}$, $\po(P',P)=w$, $P',P$ are in good position,
$\g\in U_{P'}A(P,P')/U_P$. (The associated sequence $P^n,P'{}^n$ is in this case 
$P^n=P,P'{}^n=P'$.) Thus, in this case, ${}^\tt Z_{J,\d}$ is a locally trivial 
fibration over $\cp_J$ with fibres isomorphic to $P/U_P$ for $P\in\cp_J$ so the lemma 
holds.

We now consider a general $\tt$. For any $r\ge 0$ let 
$$\tt_r=(J_n,w_n)_{n\ge r}\in\ct(J_r,\d).$$
By 3.12(b) we have a sequence of affine space bundles
$${}^\tt Z_{J,\d}@>>>{}^{\tt_1}Z_{J_1,\d}@>>>{}^{\tt_2}Z_{J_2,\d}@>>>\do\tag a$$
where for $r\gg 0$, ${}^{\tt_r}Z_{J_r,\d}$ is as in the first part of the proof. The
lemma follows.

\subhead 3.14\endsubhead
In the setup of 3.13, the maps in 3.13(a) induces bijections on the sets of $G$-orbits
(see 3.12(c)). Thus we obtain a canonical bijection between the set of $G$-orbits on 
${}^\tt Z_{J,\d}$ and the set of $G$-orbits on ${}^{\tt_r}Z_{J_r,\d}$ where $r$ is
chosen large enough so that $J_r=J_{r+1}=\do$, and $w_r=w_{r+1}=\do=w$. This last set 
of orbits is canonically the set of $(P\cap P')$-orbits on $U_{P'}\bsl A(P,P')/U_P$ 
where $P\in\cp_{J_r},P'\in\cp_{\d(J_r)}$ and $\po(P',P)=w$ (good position). Let $L_\tt$
be a common Levi of $P,P'$. Then $C_\tt=\{g\in G^1;{}^gL_\tt=L_\tt,\po({}^gP,P)=w\}$ is
a connected component of $N_{\hG}(L_\tt)$. Under the obvious bijection 
$C_\tt@>\si>>U_{P'}\bsl A(P,P)/U_P$, the conjugation action of $L_\tt$ on $C_\tt$ 
corresponds to the conjugation action of $P\cap P'$ on $U_{P'}\bsl A(P,P')/U_P$. Thus 
we obtain a canonical bijection between the set of $G$-orbits on ${}^\tt Z_{J,\d}$ and 
the set of $L_\tt$-conjugacy classes in $C_\tt$ (a connected component of an algebraic 
group with identity component $L_\tt$). Putting together these bijections we obtain a 
bijection 
$$G\bsl Z_{J,\d}\lra\sqc_{\tt\in\ct(J,\d)}L_\tt\bsl C_\tt$$
where $G\bsl Z_{J,\d}$ is the set of $G$-orbits on $Z_{J,\d}$ and $L_\tt\bsl C_\tt$ is 
the set of $L_\tt$-orbits on $C_\tt$ (for the conjugation action).

\head 4. Parabolic character sheaves on $Z_{J,\d}$ \endhead
\subhead 4.1\endsubhead   
Our notation for perverse sheaves follows \cite{\BBD}. For an algebraic variety $X$ 
over $\kk$ we write $\cd(X)$ instead of $\cd^b_c(X,\bbq)$; $l$ is a fixed prime number 
invertible in $\kk$. If $f:X@>>>Y$ is a smooth morphism with connected fibres of
dimension $d$, and $K$ is a perverse sheaf on $Y$, we set $\tf(K)=f^*(K)[d]$, a 
perverse sheaf on $X$. If $K\in\cd(X)$ and $A$ is a simple perverse sheaf on $X$ we 
write $A\dsv K$ instead of "$A$ is a composition factor of ${}^pH^i(K)$ for some 
$i\in\ZZ$".

We preserve the setup of 3.1. Let $J\sub I$. Let $B^*$ be a Borel of $G$ and let $T$ be
a maximal torus of $B^*$. Let $N^1=\{n\in G^1;nTn\i=T\}$. The map 
$N^1@>>>W,n\m\po(B^*,{}^nB^*)$ induces a bijection $T\bsl N^1=N^1/T@>\si>>W$; let 
$N^1_x$ be the $T$-coset in $N^1$ corresponding to $x\in W$. We choose $\dx\in N^1_x$. 
For $x\in W$ we define a morphism of algebraic varieties 
$\a_x:U_{B^*}N^1_xU_{B^*}@>>>N^1_x$ by $\a_x(unu')=n$ for $u,u'\in U_{B^*}$,
$n\in N^1_x$. Consider the diagram
$$N^1_x@<\ph<<\hY_x@>\r>>Y_x@>\p>>Z_{J,\d}$$
where   
$$Y_x=\{(B,B',g)\in\cz_{\em,\d};\po(B,B')=x\},$$
$$\hY_x=\{(hU_{B^*},g)\in G/U_{B^*}\T G^1;h\i gh\in U_{B^*}N^1_xU_{B^*}\},$$
$$\r(hU_{B^*},g)=({}^hB^*,{}^{gh}B^*,g),\qua \ph(hU_{B^*},g)=\a_x(h\i gh),$$
$$\p(B,B',g)=(P,P',U_{P'}gU_P),B\sub P,B'\sub P'.$$
Let $\cs(T)$ be the set of isomorphism classes of $\bbq$-local systems $\cl$ of rank 
$1$ on $T$ such that $\cl^{\ot m}\cong\bbq$ for some integer $m\ge 1$ invertible in 
$\kk$. For $\cl\in\cs(T)$ let $W^1_\cl=\{x\in W;\Ad(\dx^*)\cl\cong\cl\}$. (We have 
$\Ad(\dx):T@>>>T$.) Let $\cl\in\cs(T),x\in W^1_\cl$. The inverse image $\cl_{\dx}$ of 
$\cl$ under $N^1_x@>>>T,n\m\dx\i n,$ is a $T$-equivariant local system on $N^1_x$ for 
the conjugation action of $T$ on $N^1_x$. Now $T$ acts on $\hY_x$ by 
$t:(x,g)\m(xt\i,g)$ and on $Y_x$, trivially. The $T$-actions are compatible with 
$\ph,\r,\p$. Since $\cl_{\dx}$ is $T$-equivariant, $\ph^*\cl_{\dx}$ is a 
$T$-equivariant local system on $\hY_x$; it is of the form $\r^*\tcl$ for a well 
defined local system $\tcl$ on $Y_x$, since $\r$ is a principal $T$-bundle. We set
$$K^\cl_x=\p_!\tcl\in\cd(Z_{J,\d}).$$
Let $\hcl$ be the simple perverse sheaf on $\cz_{\em,\d}$ whose support is the closure 
of $Y_x$ and whose restriction to $Y_x$ is a shift of $\tcl$. Let 
$\bK^\cl_x=\hat\p_!\hcl$ where $\hat\p:\cz_{\em,\d}@>>>Z_{J,\d}$ is given by the same
formula as $\p$.

Occasionally we shall write $\sp$ for a local system of type $\tcl$ (as above) on $Y_x$
or on a subvariety of $Y_x$.

\subhead 4.2\endsubhead
Let $\xx=(x_1,x_2,\do,x_r)$ be a sequence in $W$, let $x=x_1x_2\do x_r$ and let
$$\align Y_\xx=&\{(B_0,B_1,\do,B_r,g)\in\cb\T\cb\T\do\T\cb\T G^1;\po(B_{i-1},B_i)=x_i,
                                               i\in[1,r],\\&B_r={}^gB_0\}.\endalign$$
Define $\p_\xx:Y_\xx@>>>Z_{J,\d}$ by $(B_0,B_1,\do,B_r,g)\m(Q,Q',g)$ where 
$Q\in\cp_J,Q'\in\cp_{\d(J)}$ are given by $B_0\sub Q,B_r\sub Q'$. Let 
$$\align\hY_\xx&=\{(h_0U_{B^*},h_1B^*,\do,h_rB^*,g)\in 
G/U_{B^*}\T G/B^*\T\do\T G/B^*\T G^1;\\&h_{i-1}\i h_i\in B^*x_iB^* \text{ for } 
i\in[1,r], h_r\i gh_0\in N_{G^1}(B^*)\}.\endalign$$ 
Define $\r_\xx:\hY_\xx@>>>Y_\xx, \ph_\xx:\hY_\xx@>>>N^1_x$ by
$$\align&
\r_\xx(h_0U_{B^*},h_1B^*,\do,h_rB^*,g)=({}^{h_0}B^*,{}^{h_1}B^*,\do,{}^{h_r}B^*,g)\\&
\ph_\xx(h_0U_{B^*},h_1B^*,\do,h_rB^*,g)=n_1n_2\do n_rn\endalign$$
where $n_i\in N_G(T)$ is given by $h_{i-1}\i h_i\in U_{B^*}n_iU_{B^*}$ and 
$n\in N^1\cap N_{G^1}(B^*)$ is given by $h_r\i gh_0\in U_{B^*}n$. Now $\ph_\xx$ is 
$T$-equivariant where $T$ acts on $\hY_\xx$ by 
$$t:(h_0U_{B^*},h_1B^*,\do,h_rB^*,g)\m(h_0t\i U_{B^*},h_1B^*,\do,h_rB^*,g).$$
Hence, if $\cl\in\cs(T)$ and $x\in W^1_\cl$ then $\ph_\xx^*(\cl_{\dx})$ is 
$T$-equivariant (see 4.1); it is of the form $\r_\xx^*\tcl$ for a well defined local 
system $\tcl$ on $Y_\xx$, since $\r_\xx$ is a principal $T$-bundle. We set
$$K^\cl_\xx=(\p_\xx)_!\tcl\in\cd(Z_{J,\d}).$$
In the case where $\xx$ reduces to a single element $x$, we have clearly
$K^\cl_\xx=K^\cl_x$. In general, $Y_\xx$ is smooth and connected. An equivalent 
statement is that
$$\align\{(h_0,h_1,\do,h_r,g)\in G^{r+1}\T G^1&;h_{i-1}\i h_i\in B^*x_iB^* (i\in[1,r]),
\\&      h_r\i gh_0B^*=B^*h_r\i gh_0\}\endalign$$
is smooth and connected. By the substitution $n=h_r\i gh_0,h_{i-1}\i h_i=y_i$, 
$i\in[1,r]$, we are reduced to the statement that
$$\{(h_0,y_1,y_2,\do,y_r,n)\in G^{r+1}\T G^1;{}^nB^*=B^*,y_i\in B^*x_iB^*\}$$
is smooth and connected, which is clear.

\subhead 4.3\endsubhead
Let $\xx=(x_1,x_2,\do,x_r)$ be a sequence in $I\sqc\{1\}$ and let $\cl\in\cs(T)$ be 
such that $x_1x_2\do x_r\in W^1_\cl$. Let
$$\align\bY_\xx=&\{(B_0,B_1,\do,B_r,g)\in\cb\T\cb\T\do\T\cb\T G^1;\\&
\po(B_{i-1},B_i)\in\{x_i,1\},i\in[1,r],B_r={}^gB_0\}.\endalign$$
Define $\bpi_\xx:\bY_\xx@>>>Z_{J,\d}$ by $(B_0,B_1,\do,B_r,g)\m(Q,Q',g)$ where 
$Q\in\cp_J,Q'\in\cp_{\d(J)}$ are given by $B_0\sub Q,B_r\sub Q'$. Now $Y_\xx$ (see 4.2)
is an open dense subset of $\bY_\xx$. By 4.2, $Y_\xx$ carries a natural local system 
$\tcl$ and the intersection cohomology complex $IC(\bY_\xx,\tcl)$ is a constructible 
sheaf $\bar\cl$ on $\bY_\xx$ (this is shown in \cite{\CS, 2.7,2.8} in the case where 
$G^1=G$; a similar proof holds in general). We set 

$\bK^\cl_\xx=(\bpi_\xx)_!\bar\cl\in\cd(Z_{J,\d})$.

\proclaim{Proposition 4.4} Let $\cl\in\cs(T)$ and let $A$ be a simple perverse sheaf on
$Z_{J,\d}$. The following conditions on $A$ are equivalent:

(i) $A\dsv K^\cl_x$ for some $x\in W^1_\cl$.

(ii) $A\dsv K^\cl_\xx$ for some sequence $\xx=(x_1,x_2,\do,x_r)$ in $W$ with
$x_1x_2\do x_r\in W^1_\cl$.

(iii) $A\dsv K^\cl_\xx$ for some sequence $\xx=(x_1,x_2,\do,x_r)$ in $I\sqc\{1\}$ with
$x_1x_2\do x_r\in W^1_\cl$.

(iv) $A\dsv\bK^\cl_\xx$ for some sequence $\xx=(x_1,x_2,\do,x_r)$ in $I\sqc\{1\}$ with
$x_1x_2\do x_r\in W^1_\cl$.

(v) $A\dsv\bK^\cl_x$ for some $x\in W^1_\cl$.
\endproclaim
In the case where $G^1=G$, the proof follows word by word that in 
\cite{\CS, 2.11-2.16}, \cite{\CSS, 12.7}. The general case can be treated in a similar
way.

\subhead 4.5\endsubhead
Let $\cc^{\cl}_{J,\d}$ be the set of (isomorphism classes) of simple perverse sheaves 
on $Z_{J,\d}$ which satisfy the equivalent conditions 4.4(i)-(v) with respect to $\cl$.
The simple perverse sheaves on $Z_{J,\d}$ which belong to $\cc^{\cl}_{J,\d}$ for some 
$\cl\in\cs(T)$ are called {\it parabolic character sheaves}; they (or their isomorphism
classes) form a set $\cc_{J,\d}$. In particular, the notion of parabolic character 
sheaf on $Z_{I,\d}=G^1$ is well defined; we thus recover the definition of character 
sheaves in \cite{\CS},\cite{\IC}.

We describe the set $\cc_{\em,\d}$. For any $x,\cl$ as in 4.1, the simple 
perverse sheaf $\hcl_x$ on $\cz_{\em,\d}$ (see 4.1) is a shift of the inverse image 
under the obvious map $\cz_{\em,\d}@>>>Z_{\em,\d}$ of a well defined simple perverse 
sheaf on $Z_{\em,\d}$. These simple perverse sheaves on $Z_{\em,\d}$ constitute 
$\cc_{\em,\d}$.

\subhead 4.6\endsubhead
Let $\tt=(J_n,w_n)_{n\ge 0}\in\ct(J,\d)$. For $r\gg 0$ we have $J_r=J_{r+1}=\do$, and 
$w_r=w_{r+1}=\do=w$. For such $r$ we define a class $\cc'_{\tt_r,\d}$ of simple 
perverse sheaves on ${}^{\tt_r}Z_{J_r,\d}$. Let $P\in\cp_{J_r},P'\in\cp_{\d(J_r)}$ and 
$\po(P',P)=w$ (good position). Let $L$ be a common Levi of $P',P$. Then 
$C=\{g\in G^1;{}^gL=L,{}^gP=P'\}$ is a connected component of $N_{\hG}(L)$. Let $X$ be 
a character sheaf on $C$ (the definition in 4.5 is applicable since $C$ is a connected 
component of an algebraic group with identity component $L$). Now $X$ is 
$L$-equivariant for the conjugation action of $L$ hence also $P\cap P'$ equivariant
($P\cap P'$ acts via its quotient $(P\cap P')/U_{P\cap P'}=L$). Hence there is a well 
defined simple perverse sheaf $X'$ on $G\T_{P\cap P'}C$ (here $P\cap P'$ acts on $G$ by
right translation) whose inverse image under $G\T C@>>>G\T_{P\cap P'}C$ is a shift of 
the inverse image of $X$ under $pr_2:G\T C@>>>C$. We may regard $X'$ as a simple 
perverse sheaf on ${}^{\tt_r}Z_{J_r,\d}$ via the isomorphism 
$$G\T_{P\cap P'}C@>\si>>{}^{\tt_r}Z_{J_r,\d},(g,c)\m({}^gP,{}^gP',gU_{P'}cU_Pg\i).$$
Now let $\vt:{}^\tt Z_{J,\d}@>>>{}^{\tt_r}Z_{J_r,\d}$ be a composition of maps in 
3.13(a); thus $\vt$ is smooth with connected fibres. Then $\tX=\ti\vt(X')$ is a simple 
perverse sheaf on ${}^\tt Z_{J,\d}$. Let $\hX$ be the simple perverse sheaf on 
$Z_{J,\d}$ whose support is the closure in $Z_{J,\d}$ of $\supp\tX$ and whose 
restriction to ${}^\tt Z_{J,\d}$ is $\tX$.

Let $\cc'_{\tt,\d}$ be the class of simple perverse sheaves on ${}^\tt Z_{J,\d}$ 
consisting of all $\tX$ as above. (It is independent of the choice of $r$.) Let 
$\cc'_{J,\d}$ be the class of simple perverse sheaves on $Z_{J,\d}$ consisting of all 
$\hX$ as above. The set of isomorphism classes of objects in $\cc'_{J,\d}$ is in 
bijection with the set of pairs $(\tt,X)$ where $\tt\in\ct(J,\d)$ and $X$ is a 
character sheaf on $C$ (as above). 

\subhead 4.7\endsubhead
We fix $b\in W$. Let $V$ be a locally closed subvariety of $Y_b$ (see 4.1). For 
$x,z\in W$ let 
$$\align X_{x,z}=&\{(\tB,B,\tB',B',g)\in\cb^4\T G^1;(\tB,\tB',g)\in V,{}^gB=B',\\&
\po(\tB,B)=x,\po(B,\tB')=z,\po(\tB',B')=\d(x)\}.\endalign$$
Define $\k:X_{x,z}@>>>V$ by $\k(\tB,B,\tB',B',g)=(\tB,\tB',g)$. Let $\cl\in\cs(T)$ be 
such that $z\d(x)\in W^1_\cl$. The inverse image of the local system $\tcl$ on 
$Y_{z,\d(x)}$ (see 4.2) under 
$$X_{x,z}@>>>Y_{z,\d(x)},\qua(\tB,B,\tB',B',g)\m(B,\tB',B',g)$$ 
is denoted again by $\tcl$.

\proclaim{Lemma 4.8} Let $m:V@>>>Z$ be a morphism of varieties. Let $A$ be a simple
perverse sheaf on $Z$ such that $A\dsv(m\k)_!\tcl$. Then there exists $\cl_1\in\cs(T)$ 
such that $b\in W^1_{\cl_1}$ and such that $A\dsv m_!\tcl_1$. (Here $\tcl_1$ is the 
local system on $Y_b$ (or its restriction to $V$) defined as in 4.1 for $b,\cl_1$ 
instead of $x,\cl$.)
\endproclaim
We argue by induction on $l(z)$. If $l(z)=0$ then $x=b$, $\k$ is an isomorphism and the
result is obvious. Assume now that $l(z)>0$. We can find $s\in W$ such that 
$l(s)=1,l(z)>l(sz)$.

Assume first that $l(xs)=l(x)+1$. Consider the isomorphism 
$$\io:X_{x,z}@>\si>>X_{xs,sz},\qua(\tB,B,\tB',B',g)\m(\tB,B_1,\tB',B'_1,g)$$ 
where $B_1\in\cb$ is given by 
$$\po(B,B_1)=s,\po(B_1,\tB')=sz\tag a$$
and $B'_1={}^gB_1$. We have $\k=\k'\io$ where $\k':X_{xs,sz}@>>>V$ is given by 
$$\k'(\tB,B_1,\tB',B'_1,g)=(\tB,\tB',g).$$
Now $\io_!\tcl$ is a local system of the same type as $\tcl$ (relative to 
$sz\d(x)\d(s)$ instead of $z\d(x)$) and $A\dsv(m\k')_!\io_!\tcl$. Using the induction 
hypothesis for $xs,sz,\io_!\tcl$ instead of $x,z,\tcl$, we get the conclusion of the 
lemma.

Assume next that $l(xs)=l(x)-1$. Then we have a partition 
$X_{x,z}=X'_{x,z}\sqc X''_{x,z}$ where $X'_{x,z}$ is the open subset defined by 
$\po(\tB,B_1)=x$ (and $B_1$ is given by (a)) and $X''_{x,z}$ is the closed subset 
defined by $\po(\tB,B_1)=xs$ (and $B_1$ is given by (a)). Let 
$j'=\k|_{X'_{x,z}},j''=\k|_{X''_{x,z}}$. By general principles, either

(b) $A\dsv(mj')_!\tcl$, or

(c) $A\dsv(mj'')_!\tcl$;
\nl
here the restriction of $\tcl$ to $X'_{x,z}$ or $X''_{x,z}$ is denoted again by $\tcl$.

Assume that (c) holds. We have $j''=\k''\io''$ where $\k'':X_{xs,sz}@>>>V$ is given by
$\k'(\tB,B_1,\tB',B'_1,g)=(\tB,\tB',g)$ and 
$$\io'':X''_{x,z}@>>>X_{xs,sz},\qua(\tB,B,\tB',B',B'_1,g)\m(\tB,B_1,\tB',B'_1,g)$$ 
(where $B_1,B'_1$ are as in $a$) is a line bundle; also $\io''_!(\tcl)$ is (up to 
shift) a local system $\tcl''$ of the same type as $\tcl$ (relative to $sz\d(x)\d(s)$ 
instead of $z\d(x)$). Hence $A\dsv(mk'')_!\tcl''$. Using the induction hypothesis for 
$xs,sz,\tcl''$ instead of $x,z,\tcl$, we get the conclusion of the lemma.

Assume now that (b) holds. We have $j'=\k'\io'$ where $\k':X_{x,sz}@>>>V$ is given by
$\k'(\tB,B_1,\tB',B'_1,g)=(\tB,\tB',g)$ and 
$$\io':X'_{x,z}@>>>X_{x,sz},\qua(\tB,B,\tB',B',B'_1,g)\m(\tB,B_1,\tB',B'_1,g)$$ 
(with $B_1,B'_1$ as in (a)) makes $X'_{x,z}$ into the complement of the zero section of
a line bundle over $X_{x,sz}$. In the case where 

(d) {\it the inverse image of $\cl$ under the coroot $\kk^*@>>>T$ corresponding to $s$
is} $\bbq$,
\nl
we have $\tcl|_{X'_{x,z}}=\io'{}^*\tcl'$ where $\tcl'$ is a local system on $X_{x,sz}$ 
of the same type as $\tcl$ (relative to $sz\d(x)$ instead of $z\d(x)$) hence we have an
exact triangle consisting of $\io'_!(\tcl)$, $\tcl'$ and a shift of $\tcl'$. Hence 
$A\dsv(mk')_!\tcl'$. Using the induction hypothesis for $x,sz,\tcl'$ instead of 
$x,z,\tcl$, we get the conclusion of the lemma. 

If (d) does not hold, the direct image $\io'_!(\tcl)$ is $0$ hence 
$${}^pH^i((mj')_!\tcl)={}^pH^i(m_!k'_!\io'_!\tcl)=0$$
and we have a contradiction. The lemma is proved.

\subhead 4.9\endsubhead
Let $\tt=(J_n,w_n)_{n\ge 0}\in\ct(J,\d)$. Let $Y_\tt$ be the set of all 
$(B,B',g)\in\cz_{\em,\d}$ such that, if $P\in\cp_J,P'\in\cp_{\d(J)}$ are given by 
$B\sub P,B'\sub P'$, then $(P,P',g)\in{}^\tt\cz_{J,\d}$. Then $(B,B',g)\m(P,P',g)$ is a
morphism $\x_\tt:Y_\tt@>>>{}^\tt\cz_{J,\d}$ and \lb $(B,B',g)\m(P,P',U_{P'}gU_P)$ is a 
morphism $\x'_\tt:Y_\tt@>>>{}^\tt Z_{J,\d}$. 

For $a\in W$ let $Y_{\tt,a}=\{(B,B',g)\in Y_\tt;\po(B,B')=a\}$. Let 
$\x_{\tt,a},\x'_{\tt,a}$ be the restrictions of $\x_\tt,\x'_\tt$ to $Y_{\tt,a}$.

\proclaim{Lemma 4.10} Let $A'$ be a simple perverse sheaf on ${}^\tt Z_{J,\d}$ with \lb
$A'\dsv(\x'_{\tt,a})_!\sp$. (Notation of 4.1.) Then 
$A'\dsv\vt{}^*(\x'_{\tt_1,b})_!\sp$ for some $b\in W$ 
($\vt:{}^\tt Z_{J,\d}@>>>{}^{\tt_1}Z_{J_1,\d}$ as in 3.12(b)).
\endproclaim
It suffices to prove the following variant of the lemma.

($*$) {\it Let $A$ be a simple perverse sheaf on ${}^\tt\cz_{J,\d}$ such that 
$A\dsv(\x_{\tt,a})_!\sp$. Then there exists $b\in W$ such that
$A\dsv\vt'{}^*(\x_{\tt_1,b})_!\sp$ where 
$\vt':{}^\tt\cz_{J,\d}@>\si>>{}^{\tt_1}\cz_{J_1,\d}$ is as in 3.12(a).}

Define $f:Y_{\tt,a}@>>>Y_{\tt_1}$ by $f(B,B',g)=({}^{g\i}(\tP'{}^B),\tP'{}^B,g)$ where 
$\x_{\tt,a}(B,B',g)=(P,P',g)$ and $\tP'=P'{}^P$. We have a partition 
$Y_{\tt_1}=\sqc_{b\in W}Y_{\tt_1,b}$. Setting $Y_{\tt,a,b}=f\i(Y_{\tt_1,b})$ we get a
partition $Y_{\tt,a}=\sqc_{b\in W}Y_{\tt,a,b}$ into locally closed subvarieties. Let 
$\x_{\tt,a,b}:Y_{\tt,a,b}@>>>{}^\tt\cz_{J,\d}$ be the restriction of $\x_{\tt,a}$. If
$A$ is as in $(*)$ then, by general principles, $A\dsv(\x_{\tt,a,b})_!\sp$ for some 
$b\in W$. We have $\vt'\x_{\tt,a}=\x_{\tt_1}f$. Hence, if 
$f_b:Y_{\tt,a,b}@>>>Y_{\tt_1,b}$ is the restriction of $f$, we have 
$\vt'\x_{\tt,a,b}=\x_{\tt_1,b}f_b$ hence $\x_{\tt,a,b}=\vt'{}\i\x_{\tt_1,b}f_b$. Thus, 
$A\dsv(\vt'{}\i)_!(\x_{\tt_1,b})_!(f_b)_!\sp$ and
$$\vt'_!A\dsv(\x_{\tt_1,b})_!(f_b)_!\sp.\tag a$$
We can write uniquely $a=a_2a_1$ where $a_1\in W_{\d(J)},a_2\in W^{\d(J)}$. We show 
that for $(\tB,\tB',g)\in Y_{\tt_1,b}$, we have
$$\align&f_b\i(\tB,\tB',g)\\&=
\{(B,B',g);B,B'\in\cb,{}^gB=B',\po(B,\tB')=a_2,\po(\tB,B)=\d\i(a_1)\}.\tag b\endalign$$
Assume first that $(B,B',g)\in f_b\i(\tB,\tB',g)$. Define $P\in\cp_J,P'\in\cp_{\d(J)}$ 
by $B\sub P,B'\sub P'$. Set $\tP'=P'{}^P$. We know that $\po(B,B')=a,\tB'=\tP'{}^B$.
We have $\po(P'{}^B,B')\in W_{\d(J)}$ (both $B',P'{}^B$ are contained in $P'$) and 
$\po(B,P'{}^B)=\po(B,P')\in W^{\d(J)}$. We have automatically $\po(P'{}^B,B')=a_1$, 
$\po(B,P'{}^B)=a_2$. It also follows that $\po({}^{g\i}(P'{}^B),B)=\d\i(a_1)$. We show 
that $P'{}^B=\tP'{}^B$. We have $P'\cap B=P'\cap P\cap B\sub P'{}^P\cap B=\tP'\cap B$. 
Also $U_{P'}\sub U_{\tP'}$ (since $\tP'\sub P'$). Hence 
$(P'\cap B)U_{P'}\sub(\tP'\cap B)U_{\tP'}$, that is $P'{}^B\sub\tP'{}^B$; since 
$P'{}^B,\tP'{}^B$ are Borels, we have $P'{}^B=\tP'{}^B=\tB'$. We see that 
$\po(B,\tB')=a_2$, $\po({}^{g\i}\tB',B)=\d\i(a_1)$. Thus, $(B,B',g)$ belongs to the 
right hand side of (b).

Conversely, assume that $(B,B',g)$ belongs to the right hand side of (b). Since 
$l(a_2a_1)=l(a_2)+l(a_1)$ we have $\po(B,B')=a_2a_1=a$. Define $P\in\cp_J$,
$P'\in\cp_{\d(J)}$ by $B\sub P,B'\sub P'$. Set $\tP'=P'{}^P$. As in the earlier part of
the proof we see that $\po(P'{}^B,B')\in W_{\d(J)}$, $\po(B,P'{}^B)\in W^{\d(J)}$ and
$P'{}^B=\tP'{}^B$. This forces $P'{}^B=\tB'$ hence $\tB'=\tP'{}^B$. Thus, 
$(B,B',g)\in f_b\i(\tB,\tB',g)$, proving (b).

We see that we may identify $Y_{\tt,a,b}$ with $X_{\d\i(a_1),a_2}$ defined as in 4.7 
relative to $V=Y_{\tt_1,b}$. Moreover, $f_b$ may be identified with $\k$ in 4.7. Now 
$\vt'_!A$ (a simple perverse sheaf on ${}^{\tt_1}\cz_{J_1,\d}$) satisfies the 
hypothesis of Lemma 4.8 with $Z={}^{\tt_1}\cz_{J_1,\d}$ and $m=\x_{\tt_1,b}$ (see (a)).
Applying Lemma 4.8, we see that there exists $\cl_1\in\cs(T)$ such that 
$b\in W^1_{\cl_1}$ and such that $\vt'_!A\dsv(\x_{\tt_1,b})_!\tcl_1$ ($\tcl_1$ as in 
4.8). Since $\vt'$ is an isomorphism, we see that 
$A\dsv\vt'{}^*(\x_{\tt_1,b})_!\tcl_1$ and $(*)$ is proved. The lemma is proved.

\proclaim{Lemma 4.11} Let $\cl\in\cs(T)$ and let $a\in W^1_\cl$. Let 
$\tt=(J_n,w_n)_{n\ge 0}\in\ct(J,\d)$. Let $\x'_{\tt,a}:Y_{\tt,a}@>>>{}^\tt Z_{J,\d}$ be
as in 4.9. Let $\tcl$ be as in 4.1 (or its restriction to $Y_{\tt,a}$). Then any 
composition factor of $\op_i{}^pH^i((\x'_{\tt,a})_!\tcl)$ belongs to $\cc'_{\tt,\d}$. 
\endproclaim
More generally we show that the lemma holds whenever $J,\tt$ are replaced by
$J_n,\tt_n$, $n\ge 0$. First we show:

(a) if the result is true for $n=1$ then it is true for $n=0$.
\nl
Let $A'$ be a composition factor of $\op_i{}^pH^i((\x'_{\tt,a})_!\tcl)$. By Lemma 4.10 
there exists $b\in W$ such that $A'\dsv\vt^*(\x'_{\tt_1,b})_!\sp$ hence 
$A'\dsv\vt^*\op_i{}^pH^i((\x'_{\tt_1,b})_!\sp)$. (Notation of 4.1.) Since $\vt$ is an 
affine space bundle, there exists a composition factor $A''$ of 
$\op_i{}^pH^i((\x'_{\tt_1,b})_!\sp)$ such that $A'=\vt^*A''$. By our hypothesis we have
$A''\in\cc'_{\tt_1,\d}$. From the definitions we have $\vt^*A''\in\cc'_{\tt,\d}$. Thus,
(a) holds.

Similarly, if the result holds for some $n\ge 1$, then it holds for $n-1$. (The proof
is the same as for $n=1$.) In this way we see that it suffices to prove the result for
$n\gg 0$. Thus, we may assume that $J_0=J_1=\do=J$, and $w_0=w_1=\do=w$. We can write 
uniquely $a=a_2a_1$ where $a_1\in W_{\d(J)},a_2\in W^{\d(J)}$. We have 
$w\in{}^{\d(J)}W$ and $a\in W_Jw\i W_{\d(J)}=w\i W_{\d(J)}$. Thus, $a_2=w\i$.

Let $P\in\cp_J,P'\in\cp_{\d(J)}$ be such that $\po(P',P)=w$. Let $L$ be a common Levi
of $P',P$. Let $C=\{c\in G^1;{}^cL=L,{}^cP=P'\}$. Let $Y'$ be the set of all  
$(\b,\b',c)$ where $\b,\b'$ are Borels of $L$ such that $\po(\b,\b')=a_1$ (position 
relative to $L$ with Weyl group $W_{\d(J)}$) and $c\in C$ satisfies ${}^c\b=\b'$. Then 
$P\cap P'$ acts on $Y'\T U_{P'}$ by $p:(\b,\b',c,u)\m({}^l\b,{}^l\b',{}^lc,{}^pu)$ 
where $l\in L,p\in lU_{P\cap P'}$. We have a commutative diagram
$$\CD
G\T_{P\cap P'}Y'@>\si>>Y_{\tt,a}\\
@VVV                @V\x'_{\tt,a}VV \\
G\T_{P\cap P'}C@>\si>>{}^\tt Z_{J,\d}   \endCD$$
where the upper horizontal map is $(g,\b,\b',c,u)\m({}^gB,{}^gB',{}^g(cu))$ with 
$B'=\b'U_{P'},B={}^{c\i}B'$, the lower horizontal map is as in 4.6, the left vertical 
map is $(g,\b,\b',c,u)\m(g,c)$. This commutative diagram shows that any composition 
factor of $\op_i{}^pH^i((\x'_{\tt,a})_!\tcl)$ is of the form $X'$ where $X$ is a 
character sheaf on $C$ (notation of 4.6); hence it is in $\cc'_{\tt,\d}$. The lemma is 
proved.

\proclaim{Lemma 4.12} For $A\in\cc_{J,\d},\tt\in\ct(J,\d)$, we set 
${}^\tt A=A|_{{}^\tt Z_{J,\d}}$. Then any composition factor of 
$\op_i{}^pH^i({}^\tt A)$ belongs to $\cc'_{\tt,\d}$.  
\endproclaim
We can find $\cl\in\cs(T)$ and $\xx=(x_1,x_2,\do,x_r)$ as in 4.3 such that 
$x_1x_2\do x_r\in W^1_\cl$ and $A\dsv\bK^\cl_\xx$ (see 4.3, 4.5). By the decomposition 
theorem \cite{\BBD} applied to $\bpi_\xx$ (see 4.3), we have 
$\bK^\cl_\xx\cong A[m]\op K'$ for some $K'\in\cd(Z_{J,\d})$ and some $m\in\ZZ$. Hence 
$\bK^\cl_\xx|_{{}^\tt Z_{J,\d}}\cong{}^\tt A[m]\op K'_1$ for some 
$K'_1\in\cd({}^\tt Z_{J,\d})$. It suffices to show that, if $A'$ is a composition 
factor of $\op_i{}^pH^i(\bK_\xx^\cl|_{{}^\tt Z_{J,\d}})$, then $A'\in\cc'_{\tt,\d}$. As
in \cite{\CS, 2.11-2.16} we see that there exists $\cl\in\cs(T),x\in W^1_\cl$ such that
$A'\dsv K_x^\cl|_{{}^\tt Z_{J,\d}}$. Using 4.11 we have $A'\in\cc'_{\tt,\d}$. The lemma
is proved.

\proclaim{Lemma 4.13} If $A\in\cc_{J,\d}$, then $A\in\cc'_{J,\d}$.
\endproclaim
Since $Z_{J,\d}=\sqc_{J\sub I}{}^\tt Z_{J,\d}$, we can find $\tt\in\ct(J,\d)$ such that
$\supp(A)\cap{}^\tt Z_{J,\d}$ is open dense in $\supp(A)$. Then
${}^\tt A=A|_{{}^\tt Z_{J,\d}}$ is a simple perverse sheaf on ${}^\tt Z_{J,\d}$ and 
${}^\tt A\in\cc'_{\tt,\d}$ (Lemma 4.12). Now $A,{}^\tt A$ are related just as $\hX,\tX$
are related in 4.6. Hence $A\in\cc'_{J,\d}$. The lemma is proved.

\proclaim{Lemma 4.14} Let $\tt=(J_n,w_n)_{n\ge 0}\in\ct(J,\d)$. Define $a\in W^{\d(J)}$
by $a\i=w_n$ for $n\gg 0$. Let $b\in W_{\d(J_\iy)}$. Let $(B,B',g)\in Y_{ab}$. Define 
$P\in\cp_J,P'\in\cp_{\d(J)}$ by $B\sub P,B'\sub P'$.

(a) We have $(P,P',g)\in{}^\tt\cz_{J,\d}$.

(b) If $b=1$, $(B,B',g)\m(P,P',g)$ is a surjective map $Y_{ab}@>>>{}^\tt\cz_{J,\d}$.
\endproclaim
We prove (a). Recall that $w_n=\min(W_{\d(J)}a\i W_{J_n})=
\min(W_{\d(J)}b\i a\i W_{J_n})$ for $n\ge 0$. In particular, 
$\po(P',P)=\min(W_{\d(J)}b\i a\i W_J)=w_0$. Define $P^n,P'{}^n$ in terms of $(P,P',g)$
as in 3.11. We have $P'{}^1=P'{}^P\in\cp_{\d(J_1)}$, $P^1\in\cp_{J_1}$. As in the proof
of Lemma 4.10 we have $ab=a_2a_1$, $\po((P'{}^1)^B,B'))=a_1\in W_{\d(J)}$, 
$a_2\in W^{\d(J)}$. Since $a\in W^{\d(J)}$, we have $a_1=b\in W_{\d(J_\iy)}$. Thus, 
$\po((P'{}^1)^B,B')\in W_{\d(J_1)}$ hence $B'\sub P'{}^1$ and $B\sub P^1$. We have
$$\po(P'{}^1,P^1)=\min(W_{\d(J_1)}\po(B',B)W_{J_1})=\min(W_{\d(J_1)}b\i a\i W_{J_1})
=w_1.$$
By the same argument for $B,B',P^1,P'{}^1,g,\tt_1$ instead of $B,B',P,P',g,\tt$, we see
that $\po(P'{}^2,P^2)=w_2$ and $B\sub P^2$, $B'\sub P'{}^2$. (We have 
$a\i\in{}^{\d(J_1)}W$ since $\d(J_1)\sub\d(J)$.) Continuing in this way, we find 
$\po(P'{}^n,P^n)=w_n$ and $B\sub P^n,B'\sub P'{}^n$ for all $n\ge 0$. In particular, 
$(P,P',g)\in{}^\tt\cz_{J,\d}$. This proves (a).

We prove (b). Assume that $(P,P',g)\in{}^\tt\cz_{J,\d}$. Define $P^n,P'{}^n$ in terms
of $(P,P',g)$ as in 3.11. By assumption, if $n\gg 0$ we have $\po(P'{}^n,P^n)=a\i$ 
(good position). Hence 
$$pr_1:\{(B,B')\in\cb\T\cb;B'\sub P'{}^n,B\sub P^n,\po(B,B')=a\}@>>>
\{B\in\cb;B\sub P^n\}$$
is an isomorphism with inverse $B\m(B,(P'{}^n)^B)$. The condition that $(B,B')$ in the 
domain of $pr_1$ satisfies $B'={}^gB$ is that $B$ is fixed by the map 
$B\m{}^{g\i}((P'{}^n)^B)$ of the flag manifold of $P^n$ into itself. This map has at
least one fixed point. Hence there exist Borels $B\sub P^n$, $B'\sub P'{}^n$ such that 
$B'={}^gB$ and $\po(B,B')=a$. We then have $B\sub P,B'\sub P'$, proving (b).

\proclaim{Lemma 4.15} Let $\tt=(J_n,w_n)_{n\ge 0}\in\ct(J,\d)$. Define
$a_2\in W^{\d(J)}$ by $w_n=a_2\i$ for $n\gg 0$. Let $a_1\in W_{\d(J_\iy)}$ and let 
$a=a_2a_1,a'=\d\i(a_1)a_2$. Let $A'$ be a simple perverse sheaf on 
${}^{\tt_1}Z_{J_1,\d}$ such that $A'\dsv(\x'_{\tt_1,a'})_!\sp$. Then $\ti\vt A'$ (a 
simple perverse sheaf on ${}^\tt Z_{J,\d}$) is a composition factor of 
$\op_i{}^pH^i((\x'_{\tt,a})_!\sp)$. 
\endproclaim
It suffices to prove the following variant of the lemma.

($*$) {\it Let $\tt=(J_n,w_n)_{n\ge 0}\in\ct(J,\d)$. Define $a_2\in W^{\d(J)}$ by 
$w_n=a_2\i$ for $n\gg 0$. Let $a_1\in W_{\d(J_\iy)}$ and let $a=a_2a_1$,
$a'=\d\i(a_1)a_2$. Let $A$ be a simple perverse sheaf on ${}^{\tt_1}\cz_{J_1,\d}$ such 
that $A\dsv(\x_{\tt_1,a'})_!\sp$. Then $\vt'{}^*A$ (a simple perverse sheaf on 
${}^\tt\cz_{J,\d}$) is a composition factor of $\op_i{}^pH^i((\x_{\tt,a})_!\sp)$.}

Define $f:Y_{\tt,a}@>>>Y_{\tt_1}$ as in the proof of Lemma 4.10. Assume that  
$f(B,B',g)=(\tB,\tB',g)$. As in the proof of 4.10 we have
$$\po(\tB,B)=\d\i(a_1),\po(B,\tB')=a_2,\po(\tB',B')=a_1;$$
By 2.6 we have $\Ad(a_2\i)J_\iy=\d(J_\iy)$ hence $a_2\i\d\i(a_1)a_2\in W_{\d(J_\iy)}$ 
and \lb
$l(a_2\i\d\i(a_1)a_2)=l(\d\i(a_1))$; since $a_2\in W^{\d(J)}\sub W^{\d(J_\iy)}$, we
have
$$l(\d\i(a_1)a_2)=
l(a_2(a_2\i\d\i(a_1)a_2))=l(a_2)+l(a_2\i\d\i(a_1)a_2)=l(a_2)+l(\d\i(a_1)).$$
Hence $\po(\tB,\tB')=\d\i(a_1)a_2=a'$. We see that $f$ defines a map
$f':Y_{\tt,a}@>>>Y_{\tt_1,a'}$. Define $f'':Y_{\tt_1,a'}@>>>Y_{\tt,a}$ by
$f(\tB,\tB',g)=(B,B',g)$ with $B\in\cb$ given by $\po(\tB,B)=\d\i(a_1),\po(B,\tB')=a_2$
and $B'={}^gB$. (From $\po(B,\tB')=a_2,\po(\tB',B')=a_1$ and $l(a_2a_1)=l(a_2)+l(a_1)$ 
we deduce that $\po(B,B')=a_2a_1$.) From the proof of 4.10 we see that $f''$ is an 
inverse to $f'$. From the definitions we have a commutative diagram
$$\CD
Y_{\tt,a}@>f'>>Y_{\tt_1,a'}\\
@V\x_{\tt,a}VV           @V\x_{\tt_1,a'}VV\\          
{}^\tt\cz_{J,\d}@>\vt'>>{}^{\tt_1}\cz_{J_1,\d}    \endCD$$
This proves ($*$) since $f'$ and $\vt'$ are isomorphisms. The lemma is proved.

\proclaim{Lemma 4.16} If $A\in\cc'_{J,\d}$, then $A\in\cc_{J,\d}$.
\endproclaim
Let $\tt=(J_n,w_n)_{n\ge 0},r,C,X,X',\tX,\hX$ be as in 4.6. We may assume that $A=\hX$.
Define $a_2\in W^{\d(J)}$ by $a_2\i=w_n$ for $n\gg 0$. Now $X\dsv(pr_3)_!\sp$ where 
$pr_3:Y'@>>>C$ (notation of 4.11).

For $r'\in[0,r]$ we denote by 
$\vt_{r'}:{}^{\tt_{r'}}Z_{J_{r'},\d}@>>>{}^{\tt_r}Z_{J_r,\d}$ the appropriate
composition of maps in 3.13(a). Consider the following statement.
 
$S_{r'}$. {\it The simple perverse sheaf $\wt{\vt_{r'}}X'$ on 
${}^{\tt_{r'}}Z_{J_{r'},\d}$ is a composition factor of  \lb
$\op_i{}^pH^i((\x'_{\tt_{r'},a_2b_{r'}})_!\sp)$ for some} $b_{r'}\in W_{\d(J_\iy)}$.

Using the commutative diagram in 4.11 with $\tt$ replaced by $\tt_r$, we see that $S_r$
holds. If $0<r'\le r$ and if $S_{r'}$ holds then $S_{r'-1}$ holds (we use Lemma 4.15 
with $\tt,\tt_1,J,J_1$ replaced by $\tt_{r'-1},\tt_{r'},J_{r'-1},J_{r'}$). Using this 
repeatedly, we see that $S_0$ holds. Thus, $\tX$ (a simple perverse sheaf on 
${}^\tt Z_{J,\d}$) is a composition factor of $\op_i{}^pH^i((\x'_{\tt,a_2b})_!\sp)$ for
some $b\in W_{\d(J_\iy)}$. Let $\p:Y_x@>>>Z_{J,\d}$ be as in 4.1 (with $x=a_2b$). By 
Lemma 4.14, $\p$ factors through a map $\p':Y_x@>>>{}^\tt Z_{J,\d}$ and 
$\tX\dsv\p'_!\sp$. Thus there exists a simple perverse sheaf on $Z_{J,\d}$ whose 
support is the closure in $Z_{J,\d}$ of $\supp\tX$, whose restriction to 
${}^\tt Z_{J,\d}$ is $\tX$ and which is a composition factor of 
$\op_i{}^pH^i(\p_!\sp)$; this is necessarily $\hX$. We see that $\hX\in\cc_{J,\d}$. The
lemma is proved.

\head 5. Central character\endhead
\subhead 5.1\endsubhead
We preserve the setup of 3.1. Let $J\sub I$. The following result is similar to
\cite{\CSS, 11.3}.

\proclaim{Lemma 5.2} Let $\cl,\cl'\in\cs(T),x\in W^1_\cl,x'\in W^1_{\cl'}$. Then
$x\in W^1_{\cl\i}$. Assume that $\cl'\not\cong\Ad(w)^*\cl$ for any $w\in W_J$. Then
$H^i_c(Z_{J,\d},K_x^{\cl\i}\ot K_{x'}^{\cl'})=0$ for all $i\in\ZZ$. Then 
\endproclaim
An equivalent statement is $\op_iH^i_c(Y_x\T_{Z_{J,\d}}Y_{x'},\wt{\cl\i}\bxt\tcl')=0$
(notation of 4.1). Here the fibre product is taken with respect to the maps
$Y_x@>>>Z_{J,\d}@<<<Y_{x'}$ as in 4.1. This fibre product is the set of all
$(B_1,B_2,B_3,B_4,g,g')$ in $\cb^4\T G^1\T G^1$ such that ${}^gB_1=B_2,{}^{g'}B_3=B_4$,
$\po(B_1,B_2)=x,\po(B_3,B_4)=x'$, $B_1,B_3$ are contained in the same parabolic $P$ of
type $J$ and $g\i g'\in U_P$. We partition this set into locally closed pieces 
$Z_w (w\in W_J)$; here $Z_w$ is defined by the condition $\po(B_1,B_3)=w$. The 
restriction of $\wt{\cl\i}\bxt\tcl'$ to a subvariety of the fibre product is denoted in
the same way. It suffices to show that $H^i_c(Z_w,\wt{\cl\i}\bxt\tcl')=0$ for all 
$w\in W_J$. For fixed $w\in W_J$, we have an obvious map $\z$ from $Z_w$ to
$$\align\{(B_1,B_2,B_3,B_4)\in\cb^4;&\po(B_1,B_2)=x,\po(B_3,B_4)=x',\po(B_1,B_3)=w,
\\&\po(B_2,B_4)=\d(w)\}.\endalign$$
Using the Leray spectral sequence of $\z$ we see that it suffices to show that $H^i_c$
of any fibre of $\z$ with coefficient in the local system above is $0$. Let
$$\Ps=\z\i({}^{g_1}B^*,{}^{g_2}B^*,{}^{g_3}B^*,{}^{g_4}B^*)$$
where $g_1,g_2,g_3,g_4\in G$. We can find $g_0\in G^1$ such that 
${}^{g_0g_1}B^*={}^{g_2}B^*,{}^{g_0g_3}B^*={}^{g_4}B^*$. We can assume that 
$g_0g_1=g_2,g_0g_3=g_4$. A point in $\Ps$ is completely determined by its 
$(g,g')$-component. Thus we may identify $\Ps$ with the set of all 
$(g,g')\in G^1\T G^1$ such that ${}^{gg_1}B^*={}^{g_0g_1}B^*$,
${}^{g'g_3}B^*={}^{g_0g_3}B^*$ and $g\i g'\in U_P$ where $P\in\cp_J$ contains
${}^{g_1}B^*$. Thus we may identify
$$\Ps=\{(g,u)\in G^1\T U_P;g_0\i g\in{}^{g_1}B^*\cap{}^{g_3}B^*\}$$
where $P$ is as above. Here $g_1,g_3\in G$ are fixed such that
$\po(B^*,{}^{g_1\i g_3}B^*)=w$. Define $\t:\Ps@>>>T$ by 
$g_1\i g_0\i gg_1\in\t(g,u)U_{B^*}$ (an affine space bundle). One checks that the local
system $\wt{\cl\i}\bxt\tcl'$ on $\Ps$ is $\t^*(\cl\i\ot\Ad(w\i)^*\cl')$. It then 
suffices to show that $H^i_c(T,\cl\i\ot\Ad(w\i)^*\cl')=0$ for all $i$. This follows 
from the fact that $\cl\i\ot\Ad(w\i)^*\cl'\in\cs(T)$ is not isomorphic to $\bbq$. The 
lemma is proved.

\subhead 5.3\endsubhead
From Lemma 5.2 we deduce as in \cite{\CSS, p.268} that there is a well defined map 
$\cc_{J,\d}@>>>\{W_J\text{-orbits in } \cs(T)\}$, $A\m(\text{$W_J$-orbit of } \cl)$ 
where $A\in\cc_{J,\d}^{\cl}$.

\head 6. The functors $f^J_{J'},e^{J'}_J$\endhead
\subhead 6.1\endsubhead
We preserve the setup of 3.1. Let $J\sub J'\sub I$. Let $Z_{J,J',\d}$ be the set of all
triples $(P,P',gU_Q)$ where $(P,P')\in\cp_J\T\cp_{\d(J)}$, $Q\in\cp_{J'}$ is given by
$P\sub Q$ and $gU_Q\in G^1/U_Q$ is such that ${}^gP=P'$ (the last condition is well 
defined since $U_Q\sub U_P$). Consider the diagram
$$Z_{J,\d}@<c<<Z_{J,J',\d}@>d>>Z_{J',\d}$$
where $c(P,P',gU_Q)=(P,P',gU_P)$ and $d(P,P',gU_Q)=(Q,Q',gU_Q)$ (with $Q$ as above and 
$Q'={}^gQ$). Define 
$$f^J_{J'}:\cd(Z_{J,\d})@>>>\cd(Z_{J',\d}),e^{J'}_J:\cd(Z_{J',\d})@>>>\cd(Z_{J,\d})$$
by $f^J_{J'}(A)=d_!c^*A,e^{J'}_J(A')=c_!d^*(A')$.

\proclaim{Lemma 6.2} For $J\sub J'\sub J''\sub I$ we have 
$f^{J'}_{J''}f^J_{J'}=f^J_{J''}$, $e^{J'}_Je^{J''}_{J'}=e^{J''}_J$.
\endproclaim
We have a diagram
$$Z_{J,\d}@<c<<Z_{J,J',\d}@>d>>Z_{J',\d}@<c'<<Z_{J',J'',\d}@>d'>>Z_{J'',\d}$$
where $c,d$ are as in 6.1 and $c',d'$ are analogous to $c,d$. We have a cartesian
diagram
$$\CD
Z_{J,J'',\d}@>d''>>Z_{J',J'',\d}\\
@Vc''VV         @Vc'VV       \\
Z_{J,J',\d}@>d>>Z_{J',\d}   \endCD$$
where $c'',d''$ are the obvious maps. Using the change of basis theorem we have
$$\align&f^{J'}_{J''}f^J_{J'}=d'_!c'{}^*d_!c^*=d'_!d''_!c''{}^* c^*=(d'd'')_!(cc'')^*=
f^J_{J''},\\& e^{J'}_Je^{J''}_{J'}=c_!d^*c'_!d'{}^*=c_!c''_!d''{}^*d'{}^*=
(cc'')_!(d'd'')^*=e^{J''}_J.\endalign$$
The lemma is proved.

\subhead 6.3\endsubhead
For $J\sub I$, let $\cd_0(Z_{J,\d})$ be the full subcategory of $\cd(Z_{J,\d})$ whose 
objects are the $A\in\cd(Z_{J,\d})$ such that any composition factor of
$\op_i{}^pH^i(A)$ belongs to $\cc_{J,\d}$.

\proclaim{Lemma 6.4} For $J\sub J'\sub I$ and $A\in\cd_0(Z_{J,\d})$, we have
$f^J_{J'}(A)\in\cd_0(Z_{J',\d})$.
\endproclaim
We may assume that $A\in\cc_{J,\d}$. Then there exists $A_0\in\cc_{\em,\d}$ such that 
$A\dsv f^\em_J(A_0)$. By the decomposition theorem \cite{\BBD}, $f^\em_J(A_0)$ is a
semisimple complex, hence some shift of $A$ is a direct summand of $f^\em_J(A_0)$. 
Hence $f^J_{J'}(A)$ is a direct summand of $f^J_{J'}f^\em_J(A_0)=f^\em_{J'}(A_0)$. 
Using the definitions and the decomposition theorem, we see that $f^\em_{J'}(A_0)$ is a
direct sum of shifts of objects in $\cc_{J',\d}$. In particular, 
$f^\em_{J'}(A_0)\in\cd_0(Z_{J',\d})$. Since $f^J_{J'}(A)$ is a direct summand of 
$f^\em_{J'}(A_0)$, we must have $f^J_{J'}(A)\in\cd_0(Z_{J',\d})$. The lemma is proved.

\proclaim{Lemma 6.5} Let $J\sub I$. If $A\in\cd_0(Z_{J,\d})$, then
$e^J_\em(A)\in\cd_0(Z_{\em,\d})$.
\endproclaim
We can assume that $A=\bK^\cl_x$ where $\cl,x$ are as in 4.1. Using the known 
relationship between $\bK^\cl_x,K^\cl_x$ (compare \cite{\CSS, 12.7}) we may assume that
$A=K^\cl_x$. For simplicity we assume also that $\cl=\bbq$. Thus, $A=K^{\bbq}_x$ where 
$x\in W$. Let 
$$\align Z&=\{(B,B',B_1,B'_1,g)\in\cb^4\T G^1;
{}^gB=B',{}^gB_1=B'_1,\\&\po(B,B')=x,\po(B,B_1)\in W_J\}.\endalign$$
Define $\ph:Z@>>>Z_{\em,\d}$ by $\ph(B,B',B_1,B'_1,g)=(B_1,B'_1,gU_{B_1})$. It suffices
to show that $\ph_!\bbq\in\cd_0(Z_{\em,\d})$. For any $z\in W_J$ let 
$$Z^z=\{(B,B',B_1,B'_1,g)\in Z;\po(B,B_1)=z\}.$$
Now $(Z^z)_{z\in W_J}$ is a partition of $Z$ into locally closed subvarieties. Let $
\ph^z:Z^z@>>>Z_{\em,\d}$ be the restriction of $\ph$. It suffices to show that 
$\ph^z_!\bbq\in\cd_0(Z_{\em,\d})$. Consider the partition 
$Z_{\em,\d}=\sqc_{w\in W}Z_{\em,\d}^w$ where
$$Z_{\em,\d}^w=\{(B_1,B'_1,gU_{B_1})\in Z_{\em,\d};\po(B_1,B'_1)=w\}.$$
It suffices to show that for any $w$, the restriction of $\ph^z_!\bbq$ to 
$Z_{\em,\d}^w$ has cohomology sheaves which are local systems with all composition 
factors isomorphic to $\bbq$. Let us fix $(B_1,B'_1)\in\cb\T\cb$ such that
$\po(B_1,B'_1)=w$. It suffices to show that the restriction of $\ph^z_!\bbq$ to 
$\{(B_1,B'_1,gU_{B_1});g\in G^1,{}^gB_1=B'_1\}$ has cohomology sheaves which are local 
systems with all composition factors isomorphic to $\bbq$. Let $T_1$ be a maximal torus
of $B_1\cap B'_1$. We can find $\z\in G^1,h\in G$ such that 
$${}^\z B_1=B'_1,{}^\z T_1=T_1,\qua {}^hT_1=T_1,\po({}^hB_1,B_1)=z.$$
Let 
$$Z'=\{(B,u,t)\in\cb\T U_{B_1}\T T_1;\po(B,{}^{\z tu}B)=x,\po(B,B_1)=z\}.$$
Consider the projection $pr'_3:Z'@>>>T_1$. We must show that $(pr'_3)_!\bbq\in\cd(T_1)$
has cohomology sheaves which are local systems with all composition factors isomorphic 
to $\bbq$. Let 
$$Z''=\{(v,u,t)\in U_{B_1}\T U_{B_1}\T T_1;\po({}^{vh}B_1,{}^{\z tuvh}B_1)=x\}.$$
Define $Z''@>>>Z'$ by $(v,u,t)\m({}^{vh}B_1,u,t)$ (an affine space bundle). Consider
the projection $pr''_3:Z''@>>>T_1$. It suffices to show that 
$(pr''_3)_!\bbq\in\cd(T_1)$ has cohomology sheaves which are local systems that are 
direct sums of copies of $\bbq$. We make the change of variable $(v,u,t)\m(v,u',t)$ 
where $u'=tuvt\i$. Since $th\in hB_1$, $Z''$ becomes $\tZ\T T_1$ where
$$\tZ=\{(v,u')\in U_{B_1}\T U_{B_1};\po({}^{vh}B_1,{}^{\z u'h}B_1)=x\}$$ 
and $pr''_3$ becomes the second projection $\tZ\T T_1@>>>T_1$. The desired conclusion 
follows. The case where $\cl\not\cong\bbq$ is treated similarly (compare with 
\cite{\GF, 2.2}.).

\proclaim{Lemma 6.6} Let $J\sub I$. If $A\in\cd(Z_{J,\d})$, then some shift of $A$ is a
direct summand of $f^\em_Je^J_\em(A)$.
\endproclaim
The argument in this proof is inspired by one in \cite{\GI, 8.5.1}. Let 
$$\cv=\{(Q,hU_Q);Q_\in\cp_J,hU_Q\in Q/U_Q\}.$$
Define $m:Z_{J,\d}\T_{\cp_J}\cv@>>>Z_{J,\d}$ by $m((P,hU_P),(P,P',gU_P))=(P,P',ghU_P)$.
For $X\in\cd(Z_{J,\d}),C\in\cd(\cv)$ we set $X\circ C=m_!(X\bxt C)\in\cd(Z_{J,\d})$. 

Let $\cv_0=\{(Q,hU_Q)\in\cv;h\in U_Q\}$ and let $j:\cv_0@>>>\cv$ be the inclusion. Let 
$\ti{\cv}$ be the set of all pairs $(B,hU_Q)$ where $B\in\cb$ and $hU_Q\in U_B/U_Q$
(with $Q\in\cp_J$ given by $B\sub Q$). Define $\p:\ti{\cv}@>>>\cv$ by 
$\p(B,hU_Q)=(Q,hU_Q)$. Clearly, the lemma is a consequence of (a),(b),(c) below.

(a) $A\circ(j_!\bbq)=A$,

(b) $A\circ(\p_!\bbq)=f^\em_Je^J_\em(A)$,

(c) $j_!\bbq[n]$ is a direct summand of $p_!\bbq$ for some $n$.
\nl
Now (a) is obvious. We prove (b). We have $f^\em_Je^J_\em(A)=d_!c^*c_!d^*(A)$ where
$$Z_{\em,\d}@<c<<Z_{\em,J,\d}@>d>>Z_{J,\d}$$
are defined by 
$$c(B,B',gU_P)=(B,B',gU_B),\qua d(B,B',gU_P)=(P,P',gU_P)$$
(with $P\in\cp_J,P'\in\cp_{\d(J)}$ given by $B\sub P,B'\sub P'$). Let $Z'$ be the set
of all quadruples $(B,B',gU_P,g'U_P)$ where $(B,B')\in\cb\T\cb,P\in\cp_J$ is given by 
$B\sub P$ and $gU_P,g'U_P\in G^1/U_P$ are such that ${}^gB=B'$ and $gU_B=g'U_B$. Let 
$b,b'$ be the projections $pr_{123}:Z'@>>>Z_{\em,J,\d}$, $pr_{124}:Z'@>>>Z_{\em,J,\d}$.
Define $a,a':Z'@>>>Z_{J,\d}$ by 
$$a(B,B',gU_P,g'U_P)=(P,P',gU_P),\qua a'(B,B',gU_P,g'U_P)=(P,P',g'U_P),$$
where $P,P'$ are as above. Then $a=db,a'=db'$. By the change of basis theorem we have 
$c^*c_!=b'_!b^*$ hence 
$$f^\em_Je^J_\em(A)=d_!b'_!b^*d^*(A)=(db')_!(db)^*(A)=a'_!a^*(A).$$
Let $Z''$ be the set of all triples $(B,gU_P,hU_P)$ where $B\in\cb,P\in\cp_J$ is given
by $B\sub P$, $gU_P\in G^1/U_P$ and $h\in U_B$. Define $\ta,\ta':Z''@>>>Z_{J,\d}$ by 
$$\ta(B,gU_P,hU_P)=(P,{}^gP,gU_P),\qua \ta'(B,gU_P,hU_P)=(P,{}^gP,ghU_P)$$
with $P$ as above. From the definitions we have $A\circ(\p_!\bbq)=\ta'_!\ta^*(A)$. The
isomorphism 
$$\io:Z''@>\si>>Z',\qua(B,gU_P,hU_P)\m(B,{}^gB,gU_P,ghU_P)$$ 
satisfies $a\io=\ta,a'\io=\ta'$. Thus, $\ta'_!\ta^*(A)=a'_!\io_!\io^*a^*(A)=a'_!a^*(A)$
and (b) is proved.         

We prove (c). Let $Q\in\cp_J$. Let 
$$\wt{Q/U_Q}=\{(B,hU_Q);B\in\cb,B\sub Q,hU_Q\in U_B/U_Q\}.$$
Let $\p':\wt{Q/U_Q}@>>>Q/U_Q$ be the second projection. It is known that $\p'$ is a
semismall map onto its image (the set of unipotent elements in $Q/U_Q$) and that 
$\p'_!\bbq$ contains as a direct summand a shift of the skyscraper sheaf of $Q/U_Q$ at
the unit element of $Q/U_Q$. Using $G$-equivariance, we see that a shift of $j_!\bbq$ 
is a direct summand of $\p_!\bbq$. This proves (c). The lemma is proved.

\proclaim{Proposition 6.7} (a) Let $J\sub I$. Let $A\in\cd(Z_{J,\d})$. We have
$A\in\cd_0(Z_{J,\d})$ if and only if $e^J_\em(A)\in\cd_0(Z_{\em,\d})$.

(b) Let $J\sub J'\sub I$. If $A'\in\cd_0(Z_{J',\d})$ then 
$e^{J'}_J(A')\in\cd_0(Z_{J,\d})$.
\endproclaim
We prove (a). If $A\in\cd_0(Z_{J,\d})$ then $e^J_\em(A)\in\cd_0(Z_{\em,\d})$ by Lemma
6.5. Conversely, assume that $e^J_\em(A)\in\cd_0(Z_{\em,\d})$. Using Lemma 6.4, we have
$f^\em_Je^J_\em(A)\in\cd_0(Z_{J,\d})$. Using this and Lemma 6.6, we see that
$A\in\cd_0(Z_{J,\d})$. 

We prove (b). Applying (a) to $A=e^{J'}_J(A')$ we see that it suffices to show that
$e^J_\em e^{J'}_J(A')\in\cd_0(Z_{\em,\d})$ or equivalently (see Lemma 6.2) that 
$e^{J'}_\em(A')\in\cd_0(Z_{\em,\d})$. But this follows from Lemma 6.5.

\proclaim{Corollary 6.8} Let $J\sub I$. Let $A$ be a character sheaf on $G^1$. Then 
$e^I_J(A)\in\cd_0(Z_{J,\d})$.
\endproclaim

\head 7. Characteristic functions\endhead
\subhead 7.1\endsubhead
We preserve the setup of 3.1. Assume that $\kk,\FF_q$ are as in 1.2 and that we are 
given an $\FF_q$-rational structure on $\hG$ with Frobenius map $F:\hG@>>>\hG$ such 
that $G^1$ is defined over $\FF_q$. Now $F$ induces on the Weyl group $W$ an 
automorphism denoted by $F:W@>>>W$; it commutes with $\d:W@>>>W$ and it carries $I$ 
onto itself. Let $J\sub I$ be such that $F(J)=J$. Then $\cp_J,\cp_{\d(J)}$ are defined 
over $\FF_q$ and $Z_{J,\d}$ is defined over $\FF_q$ with corresponding Frobenius map 
$F:Z_{J,\d}@>>>Z_{J,\d}$ given by $F(P,P',gU_P)=(F(P),F(P'),F(g)U_{F(P)})$. The natural
$G$-action on $Z_{J,\d}$ (see 3.11) restricts to a $G(\FF_q)$-action on 
$Z_{J,\d}(\FF_q)$. Let $\ce_{J,\d}$ be the vector space of all functions 
$Z_{J,\d}(\FF_q)@>>>\bbq$ that are constant on the orbits of $G(\FF_q)$. 

\subhead 7.2\endsubhead
If $X$ is an algebraic variety over $\kk$ with a fixed $\FF_q$-structure and with
Frobenius map $F:X@>>>X$ and if we are given $K\in\cd(X)$ together with an isomorphism
$\ph:F^*K@>>>K$, we denote by $\c_{K,\ph}:X(\FF_q)@>>>\bbq$ the corresponding 
characteristic function; for $x\in X(\FF_q)$, $\c_{K,\ph}(x)$ is the alternating sum 
over $i\in\ZZ$ of the trace of the map induced by $\phi$ on the stalk at $x$ of the 
$i$-th cohomology sheaf of $K$.

\subhead 7.3\endsubhead
In the setting of 7.1, let $\cc_{J,d}^F=\{A\in\cc_{J,\d};F^*A\cong A\}$. For any 
$A\in\cc_{J,\d}^F$ we choose an isomorphism $\ph_A:F^*A@>>>A$. (It is unique up to 
multiplication by an element in $\bbq^*$.) Then $\c_{A,\ph_A}:Z_{J,\d}^F@>>>\bbq$ is 
constant on the orbits of $G(\FF_q)$ since $A$ is $G$-equivariant.

Define $F:\ct(J,\d)@>\si>>\ct(J,\d)$ by 
$F((J_n,w_n)_{n\ge 0})=(F(J_n),F(w_n))_{n\ge 0}$.
Clearly, $F:Z_{J,\d}@>>>Z_{J,\d}$ carries ${}^\tt Z_{J,\d}$ onto ${}^{F(\tt)}Z_{J,\d}$.
In particular, ${}^\tt Z_{J,\d}$ is $F$-stable if and only if $\tt\in\ct(J,\d)^F$ (that
is, $F(\tt)=\tt$). For any $\tt\in\ct(J,\d)^F$, let $\ce_{\tt,\d}$ be the vector space 
of all functions ${}^\tt Z_{J,\d}(\FF_q)@>>>\bbq$ that are constant on the orbits of 
$G(\FF_q)$. We may identify $\ce_{\tt,\d}$ with a subspace $\ti{\ce}_{\tt,\d}$ of 
$\ce_{J,\d}$ by associating to $f\in\ce_{\tt,\d}$ the function $\tf\in\ce_{J,\d}$ whose
restriction to $\ce_{\tt,\d}$ is $f$ and which is $0$ on 
$Z_{J,\d}(\FF_q)-{}^\tt Z_{J,\d}(\FF_q)$. Clearly, 
$$\ce_{J,\d}=\op_{\tt\in\ct(J,\d)^F}\ti{\ce}(\tt,\d).\tag a$$
For $\tt\in\ct(J,\d)^F$ let $\cc'_{\tt,\d}{}^F=\{A'\in\cc'_{\tt,\d};F^*A'\cong A'\}$. 
For any $A'\in\cc'_{\tt,\d}{}^F$ we choose an isomorphism $\ph_{A'}:F^*A'@>>>A'$. (It 
is unique up to multiplication by an element in $\bbq^*$.) Then 
$\c_{A',\ph_{A'}}:{}^\tt Z_{J,\d}^F@>>>\bbq$ is constant on the orbits of $G(\FF_q)$. 
Thus $\c_{A',\ph_{A'}}\in\ce(\tt,\d)$ and 
$\wt{\c_{A',\ph_{A'}}}\in\ti{\ce}(\tt,\d)\sub\ce(J,\d)$.

\proclaim{Lemma 7.4} Let $A\in\cc_{J,\d}^F$. Let $d=\dim\supp(A)$. For any 
$\tt\in\ct(J,\d)$ let ${}^\tt A=A|_{{}^\tt Z_{J,\d}}$. Define $\tt^0\in\ct(J,\d)$,
$A^0\in\cc'_{\tt,\d}$ by ${}^\tt A=A^0$, ${}^\tt A=A^0$ (see Lemma 4.13). We have 
$\tt^0\in\ct(J,\d)^F,A^0\in\cc'_{\tt,\d}{}^F$ and
$$\c_{A,\ph_A}=\l\wt{\c_{A^0,\ph_{A^0}}}+\sum_{\tt\in\ct(J,\d)^F,\tt\ne\tt^0}
\sum_{A'\in\cc'_{\tt,\d}{}^F,\dim\supp(A')<d}\l_{\tt,A'}\wt{\c_{A',\ph_{A'}}}$$
where $\l_{\tt,A'}\in\bbq,\l\in\bbq^*$.
\endproclaim
For any $\tt\in\ct(J,\d)^F$, $\ph_A$ induces an isomorphism $F^*({}^\tt A)@>>>{}^\tt A$
denoted again by $\ph_A$ and we have
$$\c_{A,\ph_A}=\sum_{\tt\in\ct(J,\d)^F}\wt{\c_{{}^\tt A,\ph_A}}.$$
For $\tt=\tt^0$ we have $\wt{\c_{{}^\tt A,\ph_A}}=\l\wt{\c_{A^0,\ph_{A^0}}}$ where
$\l\in\bbq^*$. It remains to show that, for any $\tt\in\ct(J,\d)^F-\{\tt^0\}$, 
$\c_{{}^\tt A,\ph_A}$ is a $\bbq$-linear combination of functions $\c_{A',\ph_{A'}}$ 
where $A'\in\cc'_{\tt,\d}{}^F,\dim\supp(A')<d$. Clearly, 
$$\c_{{}^\tt A,\ph_A}=\sum_i(-1)^i\c_{{}^pH^i({}^\tt A),\ph_A}$$
where the isomorphism $F^*({}^pH^i({}^\tt A))@>>>{}^pH^i({}^\tt A)$ induced by $\ph_A$ 
is denoted again by $\ph_A$. It then suffices to show that 

(a) $\c_{{}^pH^i({}^\tt A),\ph_A}$ is a $\bbq$-linear combination of functions 
$\c_{A',\ph_{A'}}$ where $A'\in\cc'_{\tt,\d}{}^F,\dim\supp(A')<d$. 
\nl
From Lemma 4.13 we see that $\dim\supp{}^\tt A<d$. Hence
$\dim\supp{}^pH^i({}^\tt A)<d$ and any composition factor $A'$ of ${}^pH^i({}^\tt A)$ 
has support of dimension $<d$. By Lemma 4.12, any such $A'$ is in $\cc'_{\tt,\d}$ and 
(a) follows. The lemma is proved.

\subhead 7.5\endsubhead
We consider the following statement.

($*$) {\it The functions $\c_{A,\ph_A}$ where $A$ runs over $\cc_{J,\d}^F$ form a 
$\bbq$-basis of $\ce_{J,\d}$.}
\nl
In the case where $J=I$, so that $Z_{J,\d}=G^1$, and assuming that the characteristic 
of $\kk$ is not too small, statement ($*$) appears without proof in \cite{\IC, 5.2.1}
and with proof (when $G^1=G$) in \cite{\CSSS, 25.2}. In the remainder of this section
we assume that $(*)$ holds when $J=I$ (for $G$ and for groups of smaller dimension).

We show that ($*$) holds. By Lemma 7.4, the functions $\c_{A,\ph_A}$ are related to the
functions $\wt{\c_{A',\ph_{A'}}}$ (with $\tt\in\ct(J,\d)^F,A'\in\cc'_{\tt,\d}{}^F$) by 
an upper triangular matrix with invertible entries on the diagonal. Hence it suffices 
to show that the functions $\wt{\c_{A',\ph_{A'}}}$ (with $\tt\in\ct(J,\d)^F$,
$A'\in\cc'_{\tt,\d}{}^F$) form a $\bbq$-basis of $\ce_{J,\d}$. Using 7.3(a), we see 
that it suffices to show that for any $\tt\in\ct(J,\d)^F$, the functions 
$\c_{A',\ph_{A'}}$ (with $A'\in\cc'_{\tt,\d}{}^F$) form a $\bbq$-basis of 
$\ce_{\tt,\d}$. Using 3.14 and the definition of $\cc'_{\tt,\d}$ (see 4.6), we see that
it suffices to prove the statement in the previous sentence when ${}^\tt Z_{J,\d}$ are 
replaced by ${}^{\tt_r}Z_{J_r,\d}$ (notation of 3.14) for $r$ large. We may assume that
$P,P',L_\tt,C_\tt$ (as in 3.14) are defined over $\FF_q$. We are reduced to the 
statement that the characteristic functions of the character sheaves on $C_\tt$ 
(relative to the $\FF_q$-structure) form a $\bbq$-basis for the space of 
$L_\tt(\FF_q)$-invariant functions $C_\tt(\FF_q)@>>>\bbq$, which is part of our 
assumption.

\enddocument